\documentclass[10pt]{amsart}
\usepackage{geometry}                
\usepackage{graphicx,subfigure}
\usepackage{amssymb,amsthm}
\usepackage[authoryear,sort&compress]{natbib}
\usepackage{overpic,rotating}
\bibpunct{(}{)}{;}{a}{,}{,}
\usepackage{multibib}
\newcites{pri}{Primary Literature}
\newcites{rev}{Books and Reviews}

\def\pder#1#2{{\frac {\partial #1} {\partial #2}}}

\renewcommand{\Re}{\ensuremath{\mathbb{R}}}
\newcommand{\bfi}{\bfseries\itshape}
\newcommand{\deriv}[2]{\ensuremath{\frac{\partial #1}{\partial #2}}}
\newcommand{\refeqn}[1]{(\ref{eqn:#1})}
\newcommand{\reffig}[1]{Figure \ref{fig:#1}}
\newcommand{\norm}[1]{\ensuremath{\left\| #1 \right\|}}

\newcommand{\bracket}[1]{\ensuremath{\left[ #1 \right]}}
\newcommand{\braces}[1]{\ensuremath{\left\{ #1 \right\}}}
\newcommand{\parenth}[1]{\ensuremath{\left( #1 \right)}}
\newcommand{\tr}[1]{\mbox{tr}\ensuremath{\negthickspace\bracket{#1}}}
\newcommand{\so}{\ensuremath{\mathfrak{so}(3)}}
\newcommand{\SO}{\ensuremath{\mathrm{SO}(3)}}
\newcommand{\SE}{\ensuremath{\mathrm{SE}(3)}}
\renewcommand{\Re}{\ensuremath{\mathbb{R}}}
\newcommand{\pair}[1]{\ensuremath{\langle #1 \rangle}}
\newcommand{\pairlr}[1]{\ensuremath{\left\langle #1 \right\rangle}}

\renewcommand{\S}{\ensuremath{\mathbb{S}}}

\newtheoremstyle{example}
  {10pt}
  {3pt}
  {}
  {}
  {\bfseries}
  {.}
  {\newline}
  {}
\theoremstyle{example}
\newtheorem{example}{Example}

\title{Discrete Control Systems}
\dedicatory{Invited article for the Springer  Encyclopedia of Complexity and System Science}
\author{Taeyoung Lee}
\address{Department of Aerospace Engineering, University of Michigan, Ann Arbor, Michigan, USA}
\author{Melvin Leok}
\address{Department of Mathematics, Purdue University, West Lafayette, Indiana, USA}
\author{N. Harris McClamroch}
\address{Department of Aerospace Engineering, University of Michigan, Ann Arbor, Michigan, USA}
\thanks{TL and ML have been supported in part by NSF Grant DMS-0504747 and DMS-0726263.
TL and NHM have been supported in part by NSF Grant ECS-0244977 and CMS-0555797.}

\begin{document}

\maketitle

\setcounter{tocdepth}{1} \tableofcontents

\section*{Glossary and Notation}
\begin{description}
\item[Discrete variational mechanics] A formulation of mechanics in discrete-time that is based on a discrete analogue of \textbf{Hamilton's principle}, which states that the system takes a trajectory for which the action integral is stationary.
\item[Geometric integrator] A numerical method for obtaining numerical solutions of differential equations that preserves geometric properties of the continuous flow, such as symplecticity, momentum preservation, and the structure of the configuration space.
\item[Lie group] A differentiable manifold with a group structure where the composition is differentiable. The corresponding \textbf{Lie algebra} is the tangent space to the Lie group based at the identity element. 
\item[Symplectic] A map is said to be symplectic if given any initial volume in phase space, the sum of the signed projected volumes onto each position-momentum subspace is invariant under the map. One consequence of symplecticity is that the map is volume-preserving as well.
\end{description}

\section{Definition of the Subject and Its' Importance}
Discrete control systems, as considered here, refer to the control theory of discrete-time Lagrangian or Hamiltonian systems. These discrete-time models are based on a discrete variational principle, and are  part of the broader field of geometric integration. Geometric integrators are numerical integration methods that preserve geometric properties of continuous systems, such as conservation of the symplectic form, momentum, and energy. They also guarantee that the discrete flow remains on the manifold on which the continuous system evolves, an important property in the case of rigid-body dynamics.

In nonlinear control, one typically relies on differential geometric and dynamical systems techniques to prove properties such as stability, controllability, and optimality. More generally, the geometric structure of such systems plays a critical role in the nonlinear analysis of the corresponding control problems. Despite the critical role of geometry and mechanics in the analysis of nonlinear control systems, nonlinear control algorithms have typically been implemented using numerical schemes that ignore the underlying geometry.

The field of discrete control system aims to address this deficiency by restricting the approximation to choice of a discrete-time model, and developing an associated control theory that does not introduce any additional approximation. In particular, this involves the construction of a control theory for discrete-time models based on geometric integrators that yields numerical implementations of nonlinear and geometric control algorithms that preserve the crucial underlying geometric structure.

\section{Introduction}
The dynamics of Lagrangian and Hamiltonian systems have unique geometric properties; the Hamiltonian flow is symplectic, the total energy is conserved in the absence of non-conservative forces, and the momentum map associated with a symmetry of the system is preserved. Many interesting dynamics evolve on a non-Euclidean space. For example, the configuration space of a spherical pendulum is the two-sphere, and the configuration space of a rigid body attitude dynamics has a Lie group structure, namely the special orthogonal group. These geometric features determine the qualitative behavior of the system, and serve as a basis for theoretical study.

Geometric numerical integrators are numerical integration algorithms that preserve  structures of the continuous dynamics such as invariants, symplecticity, and the configuration manifold (see \citepri{HaLuWa2006}). The exact geometric properties of the discrete flow not only generate improved qualitative behavior, but also provide accurate and efficient numerical techniques.  In this article, we view a geometric integrator as an intrinsically discrete dynamical system, instead of concentrating on the numerical approximation of a continuous trajectory.

Numerical integration methods that preserve the simplecticity of a Hamiltonian system have been studied (see \citepri{San.AN92,LeRe2004}). Coefficients of a Runge-Kutta method are carefully to chosen to satisfy a simplecticity criterion and order conditions to obtain a symplectic Runge-Kutta method. However, it can be difficult to construct such integrators, and it is not guaranteed that other invariants of the system, such as a momentum map, are preserved. Alternatively, variational integrators are constructed by discretizing Hamilton's principle, rather than discretizing the continuous Euler-Lagrange equation (see \citepri{MosVes.CMP91,MaWe2001}). The resulting integrators have the desirable property that they are symplectic and momentum preserving, and they exhibit good energy behavior for exponentially long times. Lie group methods are numerical integrators that preserve the Lie group structure of the configuration space (see \citepri{IsMuNoZa2000}). Recently, these two approaches have been unified to obtain Lie group variational integrators that preserve the geometric properties of the dynamics as well as the Lie group
structure of the configuration space without the use of local charts, reprojection, or constraints (see \citepri{MaPeSh1999,Leo.Phd04,CMA07}).

Optimal control problems involve finding a control input  such that a certain optimality objective is achieved under prescribed constraints. An optimal control problem that minimizes a performance index is described by a set of differential equations, which can be derived using Pontryagin's maximum principle. Discrete optimal control problems involve finding a control input for a discrete dynamic system such that an optimality objective is achieved with prescribed constraints. Optimality conditions are derived from the discrete equations of motion, described by a set of discrete equations. This approach is in contrast to traditional techniques where a discretization appears at the last stage to solve the optimality condition numerically. Discrete mechanics and optimal control approaches determine optimal control inputs and trajectories more accurately with less computational load (see \citepri{JuMaOb2005}). Combined with an indirect optimization technique, they are substantially more efficient (see \citepri{HuLeSaBl2006,CDC06.opt,ACC06}).

The geometric approach to mechanics can provide as the theoretical basis of innovative control methodologies in geometric control theory. For example, these techniques allow the attitude of satellites to be controlled using changes in its shape, as opposed to chemical propulsion.
While the geometric structure of mechanical systems plays a critical role in the construction of geometric control algorithms, these algorithms have typically been implemented using numerical schemes that ignore the underlying geometry. By applying geometric control algorithms to discrete mechanics that preserve geometric properties, we obtain exact numerical implementation of the geometric control theory. In particular, the method of controlled Lagrangian systems is based on the idea of adopting a feedback control to realize a modification of either the potential energy or the kinetic energy, referred to as potential shaping or kinetic shaping, respectively. These ideas are applied to
construct a real-time digital feedback controller that stabilizes the inverted equilibrium of the cart-pendulum (see \citepri{BloLeoMar.CDC05,BloLeoMar.CDC06}).

In this article, we will survey discrete Lagrangian and Hamiltonian mechanics, and their applications to optimal control and feedback control theory.

\section{Discrete Lagrangian and Hamiltonian Mechanics}\label{sec:dm}

Mechanics studies the dynamics of physical bodies acting under forces and potential fields. In Lagrangian mechanics, the trajectory of the object is derived by finding the path that minimizes the integral of a Lagrangian over time, called the action integral. In many classical problems, the Lagrangian is chosen as the difference between kinetic energy and potential energy. The Legendre transformation provides an alternative description of mechanical systems, referred to as Hamiltonian mechanics.

\begin{figure}
\setlength{\unitlength}{1.75em}\centering\small
\begin{picture}(22.3,14)(0,-14)
\put(0.0,-2.0){\framebox(5.3,2.0)[c]
{\shortstack[c]{Configuration Space\\$(q,\dot{q})\in TQ$}}}
\put(2.65,-2.0){\vector(0,-1){1.0000}}
\put(0.0,-5.0){\framebox(5.3,2.0)[c]
{\shortstack[c]{Lagrangian\\$L(q,\dot{q})$}}} \put(2.65,-5.0){\vector(0,-1){1.0000}}
\put(2.65,-5.5){\line(1,0){5.5}}
\put(8.15,-5.5){\vector(0,-1){3.5}}
\put(0.0,-8.0){\framebox(5.3,2.0)[c]
{\shortstack[c]{Action Integral\\$\mathfrak{G}=\int_{t_0}^{t_f} L(q,\dot{q})\, dt$}}}
\put(2.65,-8.0){\vector(0,-1){1.0000}}
\put(0.0,-11.0){\framebox(5.3,2.0)[c]
{\shortstack[c]{Variations\\$\delta\mathfrak{G}=\frac{d}{d\epsilon}\mathfrak{G}^\epsilon=0$}}}
\put(2.65,-11.0){\vector(0,-1){1.0000}}
\put(0.0,-14.0){\framebox(5.3,2.0)[c]
{\shortstack[c]{Euler--Lagrange Eqn.\\$\frac{d}{dt}D_2 L -D_1 L =0$}}}
\put(5.5,-11.0){\framebox(5.3,2.0)[c]
{\shortstack[c]{Legendre transform.\\$p=\deriv{}{\dot q}L(q,\dot{q})$}}}
\put(8.15,-11.0){\vector(0,-1){1.0000}}
\put(5.5,-14.0){\framebox(5.3,2.0)[c]
{\shortstack[c]{Hamilton's Eqn.\\$\dot{q}=H_p,\,\dot{p}=-H_q$}}}
\put(11.5,-2.0){\framebox(5.3,2.0)[c]
{\shortstack[c]{Configuration Space\\$(q_k, q_{k+1})\in Q\times Q$}}}
\put(14.15,-2.0){\vector(0,-1){1.0000}}
\put(11.5,-5.0){\framebox(5.3,2.0)[c]
{\shortstack[c]{Discrete Lagrangian\\$L_{d}(q_k,q_{k+1})$}}}
\put(14.15,-5.0){\vector(0,-1){1.0000}}
%
\put(14.15,-5.5){\line(1,0){5.5}}
\put(19.65,-5.5){\vector(0,-1){3.5}}
\put(11.5,-8.0){\framebox(5.3,2.0)[c]
{\shortstack[c]{Action Sum\\$\mathfrak{G}_d=\sum L_{d}(q_k,q_{k+1})$}}}
\put(14.15,-8.0){\vector(0,-1){1.0000}}
\put(11.5,-11.0){\framebox(5.3,2.0)[c]
{\shortstack[c]{Variation\\$\delta\mathfrak{G}_d=\frac{d}{d\epsilon}\mathfrak{G}_d^\epsilon=0$}}}
\put(14.15,-11.0){\vector(0,-1){1.0000}}
\put(11.5,-14.0){\framebox(5.3,2.0)[c]
{\shortstack[c]{Dis. E-L Eqn.
\\$D_2 L_{d_{k-1}}\!+\!D_1 L_{d_k}=0$}}}
\put(17,-11.0){\framebox(5.3,2.0)[c]
{\shortstack[c]{Legendre transform.\\$p_k=\deriv{}{\dot q}L(q,\dot{q})\big|_k$}}}
\put(19.65,-11.0){\vector(0,-1){1.0000}}
\put(17,-14.0){\framebox(5.3,2.0)[c]
{\shortstack[c]{Dis. Hamilton's Eqn.
\\\scriptsize{$p_k=-D_1 L_{d_k},$}\\
\scriptsize{$p_{k+1}=D_2 L_{d_k}$}}}}
\end{picture}
\caption{Procedures to derive continuous and discrete equations of motion}\label{fig:el}
\end{figure}
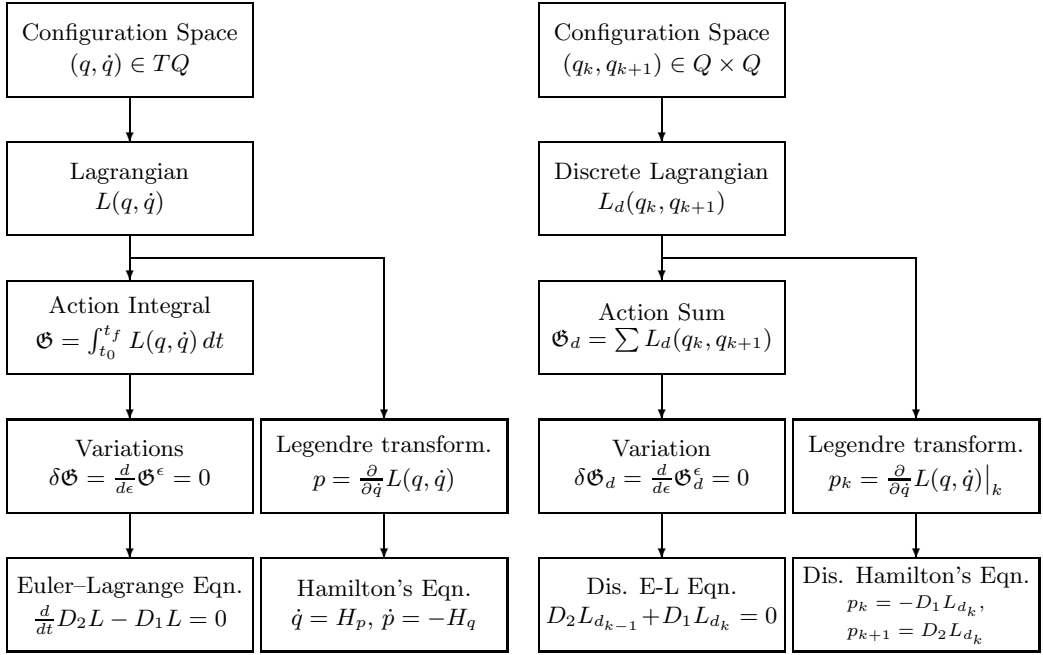

Discrete Lagrangian and Hamiltonian mechanics has been developed by reformulating the theorems and the procedures of Lagrangian and Hamiltonian mechanics in a discrete time setting (see, for example, \citepri{MaWe2001}). Therefore, discrete mechanics has a parallel structure with the mechanics described in continuous time, as summarized in \reffig{el} for Lagrangian mechanics. In this section, we describe  discrete Lagrangian mechanics in more detail, and we derive discrete Euler-Lagrange equations for several mechanical systems.

Consider a mechanical system on a configuration space $Q$, which is the space of possible positions. The Lagrangian depends on the position and velocity, which are elements of the tangent bundle to $Q$, denoted $TQ$. Let $L:TQ\rightarrow\Re$ be the Lagrangian of the system. The discrete Lagrangian, $L_d:Q\times Q\rightarrow \mathbb{R}$ is an approximation to the exact discrete Lagrangian,
\begin{align}
L_d^{\operatorname{exact}}(q_0,q_1)=\int_0^h L(q_{01}(t),\dot q_{01}(t)) dt,\label{eqn:Ld}
\end{align}
where $q_{01}(0)=q_0,$ $q_{01}(h)=q_1,$ and $q_{01}(t)$ satisfies the
Euler--Lagrange equation in the time interval $(0,h)$. A discrete action sum  $\mathfrak{G}_d:Q^{N+1}\rightarrow \mathbb{R}$, analogous to the action integral, is given by
\begin{align}
\mathfrak{G}_d(q_0,q_1,\ldots,q_N) = \sum_{k=0}^{N-1} L_d(q_k, q_{k+1}). \label{eqn:Gd}
\end{align}
The discrete Hamilton's principle states that
\[ \delta \mathfrak{G}_d = 0\]
for any $\delta q_k$, which yields the {\bfi discrete Euler--Lagrange (DEL)} equation,
\begin{align}
D_2 L_d(q_{k-1},q_k)+D_1 L_d(q_k,q_{k+1})=0.\label{eqn:DEL}
\end{align}
This yields a discrete Lagrangian flow map $(q_{k-1},q_k)\mapsto(q_k,q_{k+1})$. The discrete Legendre transformation, which from a pair of positions $(q_0, q_1)$ gives a position-momentum pair $(q_0,p_0)=(q_0,-D_1 L_d(q_0,q_1)$ provides a discrete Hamiltonian flow map in terms of momenta.

The discrete equations of motion, referred to as variational integrators, inherit the geometric properties of the continuous system; they are symplectic, and they preserve any momentum maps associated with symmetries as the discrete Noether's theorem is satisfied. They exhibit good total energy behavior for exponentially long time periods.

Many interesting Lagrangian and Hamiltonian systems, such as rigid bodies evolve on a Lie group. Lie group variational integrators preserve the nonlinear structure of the Lie group configurations as well as geometric properties of the continuous dynamics (see \citepri{MaPeSh1999} and \citepri{Leo.Phd04}). The basic idea for all Lie group methods is to express the update map for the group elements in terms of the group operation,
\begin{align}
    g_1 = g_0 f_0,\label{eqn:g1}
\end{align}
where $g_0,g_1\in G$ are configuration variables in a Lie group $G$, and $f_0 \in G$ is the discrete update represented by a right group operation on $g_0$. Since the group element is updated by a group operation, the group structure is preserved automatically without need of parameterizations, constraints, or re-projection. In the Lie group variational integrator, the expression for the flow map is obtained from the discrete variational principle on a Lie group, the same procedure presented in \reffig{el}. But, the infinitesimal variation of a Lie group element must be carefully expressed to respect the structure of the Lie group. For example, it can be expressed in terms of the exponential map as
\begin{align*}
    \delta g = \frac{d}{d\epsilon}\bigg|_{\epsilon=0} g \exp {\epsilon\eta} = g\eta,
\end{align*}
for a Lie algebra element $\eta\in\mathfrak{g}$. This approach has been applied to the rotation group $\mathrm{SO}(3)$ and to the special Euclidean group $\mathrm{SE}(3)$ for dynamics of rigid bodies (see \citepri{CCA05}, \citepri{CMA07}, and \citepri{CMDA07}). Generalizations to arbitrary Lie groups gives the generalized {\bfi discrete Euler--Poincar\'{e} (DEP)} equation,
\begin{align}
    T_e^* L_{f_{0}}\cdot D_2 L_d (g_{0},f_{0}) - \mathrm{Ad}^*_{f_{0}}\cdot(T_e^* L_{f_{1}}\cdot D_2 L_d (g_{1},f_{1})) + T_e^* L_{g_{1}}\cdot D_1 L_d (g_{1},f_{1})=0,\label{eqn:DEP}
\end{align}
for a discrete Lagrangian on a Lie group, $L_d:G\times G\rightarrow \Re$. Here $L_f:G\rightarrow G$ denotes the left translation map given by $L_f g = fg$ for $f,g\in G$, $T_g L_f:T_g G\rightarrow T_{fg} G$ is the tangential map for the left translation, and $\mathrm{Ad}_g:\mathfrak{g}\rightarrow\mathfrak{g}$ is the adjoint map. A dual map is denoted by a superscript ${}^*$ (see \citepri{MarRat.BK99} for detailed definitions).

We illustrate the properties of discrete mechanics using several mechanical systems, namely a mass-spring system, a planar pendulum, a spherical pendulum, and a rigid body.

\begin{example}[\textbf{Mass-spring System}]
Consider a mass-spring system, defined by a rigid body that moves along a straight frictionless slot, and is attached to a linear spring.

\textit{Continuous equation of motion:}
The configuration space is $Q=\Re$, and the Lagrangian $L:\Re\times\Re\rightarrow \Re$ is given by
\begin{align}
    L(q,\dot q) =\frac{1}{2}m\dot q^2 -\frac{1}{2}\kappa q^2,\label{eqn:Lms}
\end{align}
where $q\in \Re$ is the displacement of the body measured from the point where the spring exerts no force. The mass of the body and the spring constant are denoted by $m,\kappa\in\Re$, respectively. The Euler-Lagrange equation yields the continuous equation of motion.
\begin{align}
    m \ddot q + \kappa q =0.\label{eqn:ELms}
\end{align}
\textit{Discrete equation of motion:} Let $h>0$ be a discrete time step, and a subscript $k$ denotes the $k$-th discrete variable at $t=kh$. The discrete Lagrangian $L_d:\Re\times\Re\rightarrow\Re$ is an approximation of the integral of the continuous Lagrangian \refeqn{Lms} along the solution of \refeqn{ELms} over a time step. Here, we choose the following discrete Lagrangian.
\begin{align}
    L_d(q_k,q_{k+1})= hL\!\parenth{\frac{q_k+q_{k+1}}{2},\frac{q_{k+1}-q_{k}}{h}} = \frac{1}{2h} m (q_{k+1}-q_k)^2 -\frac{h\kappa}{8} (q_k+q_{k+1})^2.\label{eqn:Ldms}
\end{align}
Direct application of the discrete Euler-Lagrange equation to this discrete Lagrangian yields the discrete equations of motion. We develop the discrete equation of motion using the discrete Hamilton's principle to illustrate the principles more explicitly. Let $\mathfrak{G}_d:\Re^{N+1}\rightarrow\Re$ be the discrete action sum defined as $\mathfrak{G}_d=\sum_{k=0}^{N-1} L_d(q_k,q_{k+1})$, which approximates the action integral. The infinitesimal variation of the action sum can be written as
\begin{align*}
    \delta \mathfrak{G}_d & = \sum_{k=0}^{N-1} \delta q_{k+1} \braces{\frac{m}{h}(q_{k+1}-q_k)-\frac{h\kappa}{4}(q_k+ q_{k+1})}
    +\delta q_{k} \braces{-\frac{m}{h}(q_{k+1}-q_k)-\frac{h\kappa}{4}(q_k+ q_{k+1})}.
\end{align*}
Since $\delta q_0=\delta q_N=0$, the summation index can be rewritten as
\begin{align*}
    \delta \mathfrak{G}_d  & = \sum_{k=1}^{N-1} \delta q_{k} \braces{-\frac{m}{h}(q_{k+1}-2 q_k+ q_{k-1})-\frac{h\kappa}{4} (q_{k+1}+2 q_{k}+q_{k-1})}.
\end{align*}
From discrete Hamilton's principle, $\delta \mathfrak{G}_d =0$ for any $\delta q_k$. Thus, the discrete equation of motion is given by
\begin{align}
    \frac{m}{h}(q_{k+1}-2 q_k+ q_{k-1})+\frac{h\kappa}{4} (q_{k+1}+2 q_{k}+q_{k-1})=0.
\end{align}
For a given $(q_{k-1},q_k)$, we solve the above equation to obtain $q_{k+1}$. This yields a discrete flow map $(q_{k-1},q_k)\mapsto(q_k,q_{k+1})$, and this process is repeated. The discrete Legendre transformation provides the discrete equation of motion in terms of the velocity as
\begin{gather}
    \parenth{1+\frac{h^2\kappa}{4m}}q_{k+1} =h\dot q_k + \parenth{1-\frac{h^2\kappa}{4m}} q_k,\label{eqn:qkpms}\\
    \dot q_{k+1}  = \dot q_k -\frac{h\kappa}{2m} q_k - \frac{h\kappa}{2m} q_{k+1}.\label{eqn:dotqkpms}
\end{gather}
For a given $(q_{k},\dot q_k)$, we compute $q_{k+1}$ and $\dot q_{k+1}$ by \refeqn{qkpms} and \refeqn{dotqkpms}, respectively. This yields a discrete flow map $(q_{k},\dot q_k)\mapsto(q_{k+1},\dot q_{k+1})$. It can be shown that this variational integrator has second-order accuracy, which follows from the fact that the discrete action sum is a second-order approximation of the action integral.

\textit{Numerical example:} We compare computational properties of the discrete equations of motion given by \refeqn{qkpms} and \refeqn{dotqkpms} with a 4(5)-th order variable step size Runge-Kutta method. We choose $m=1\,\mathrm{kg}$, $\kappa=1\,\mathrm{kg/s^2}$ so that the natural frequency is $1\,\mathrm{rad/s}$. The initial conditions are $q_0=\sqrt{2}\,\mathrm{m}$, $\dot q_0=0$, and the total energy is $E=1\,\mathrm{Nm}$. The simulation time is $200\pi\,\mathrm{sec}$, and the step-size $h=0.035$ of the discrete equations of motion is chosen such that the CPU times are the same for both methods. \reffig{Ems} shows the computed total energy. The variational integrator preserves the total energy well. The mean variation is $2.7327\times 10^{-13}\,\mathrm{Nm}$. But, there is a notable dissipation of the computed total energy for the Runge-Kutta method

\begin{figure}
\subfigure[Mass-spring system]{
\includegraphics[width=0.35\textwidth]{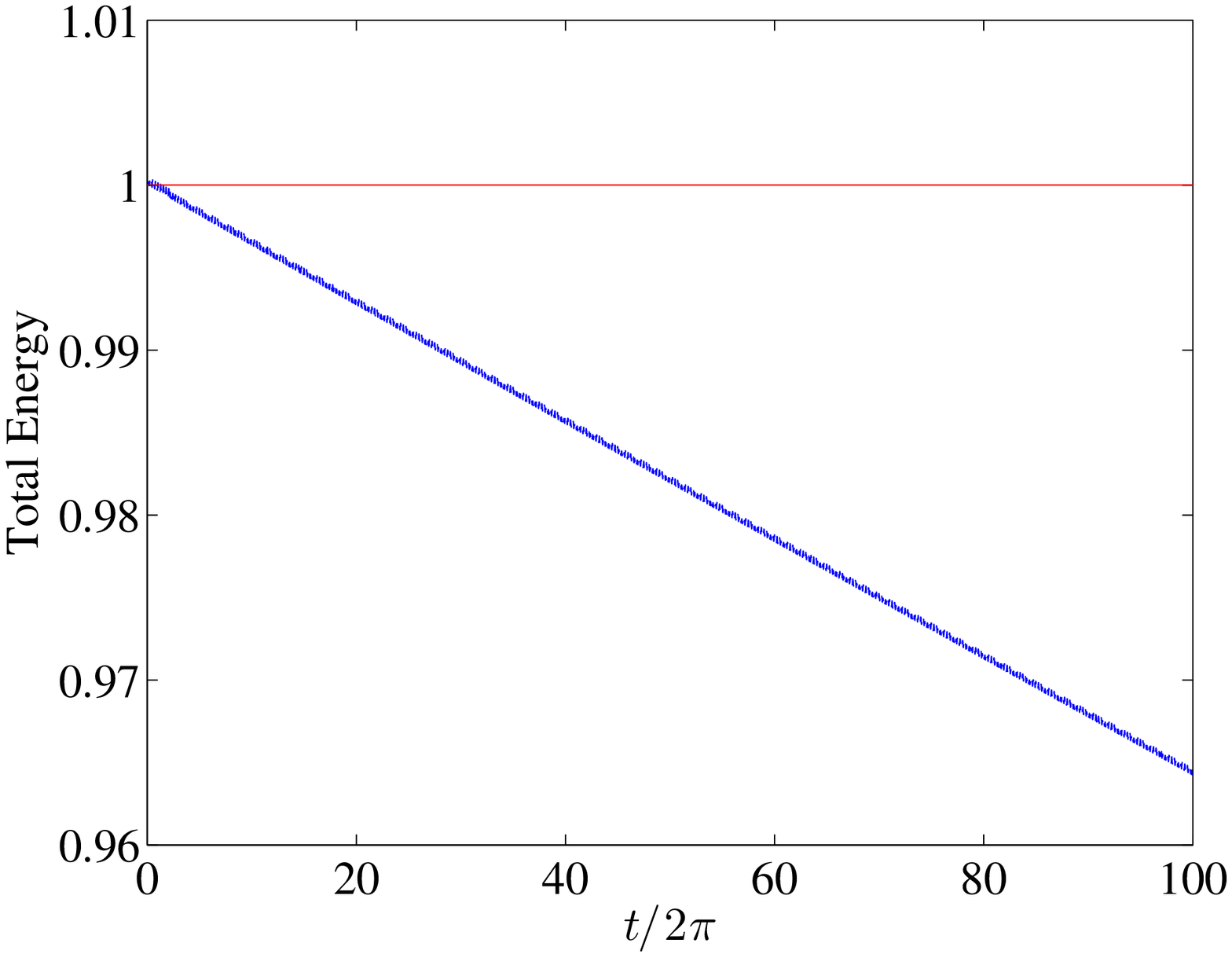}\label{fig:Ems}}
\hspace*{0.05\textwidth}
\subfigure[Planar pendulum]{
\includegraphics[width=0.36\textwidth]{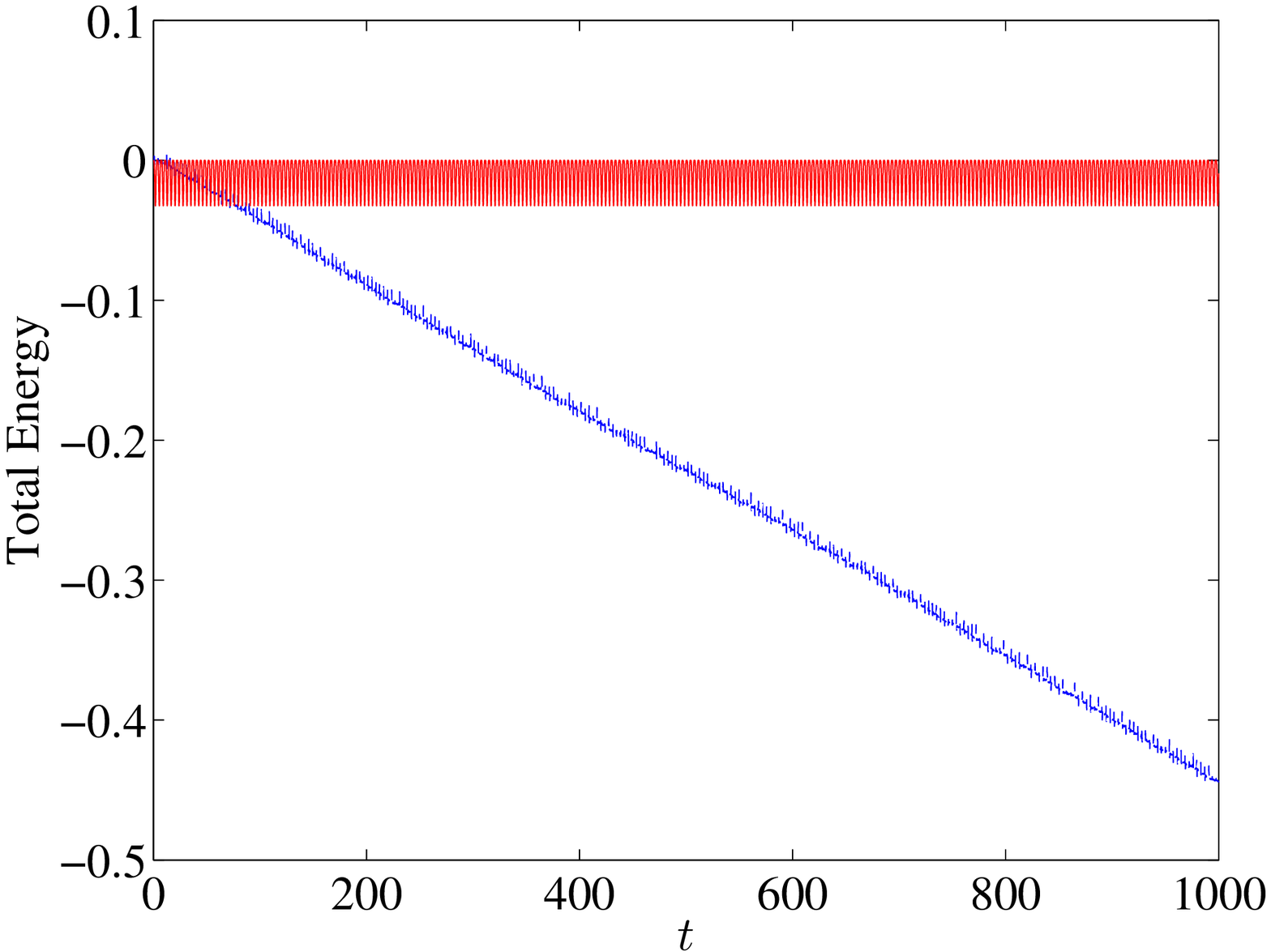}\label{fig:Epp}}
\caption{Computed total energy (RK45: blue, dotted, VI: red, solid)}
\end{figure}

\end{example}

\begin{example}[\textbf{Planar Pendulum}]\label{ex:pp}
A planar pendulum is a mass particle connected to a frictionless, one degree-of-freedom pivot by a rigid massless link under a uniform gravitational potential. The configuration space is the one-sphere $\S^1=\braces{q\in\Re^2\,|\,\norm{q}=1}$. While it is common to parameterize the one-sphere by an angle, we develop parameter-free equations of motion in the special orthogonal group $\mathrm{SO}(2)$, which is a group of $2\times 2$ orthogonal matrices with determinant 1, i.e. $\mathrm{SO}(2)=\braces{R\in\Re^{2\times 2}\,|\,R^TR=I_{2\times 2},\mathrm{det}[R]=1}$. $\mathrm{SO}(2)$ is diffeomorphic to the one-sphere. It is also possible to develop global equations of motion on the one-sphere directly, as shown in the next example, but here we focus on the special orthogonal group in order to illustrate the key steps to develop a Lie group variational integrator.

We first exploit the basic structures of the Lie group $\mathrm{SO}(2)$. Define a hat map $\hat{\cdot}$, which maps a scalar $\Omega$ to a $2\times 2$ skew-symmetric matrix $\hat\Omega$ as
\begin{align*}
    \hat\Omega=\begin{bmatrix} 0 & -\Omega \\ \Omega & 0 \end{bmatrix}.
\end{align*}
The $2\times 2$ skew-symmetric matrices are the Lie algebra $\mathfrak{so}(2)$.
Using the hat map, we identify $\mathfrak{so}(2)$ with $\Re$. An inner product on $\mathfrak{so}(2)$ can be induced from the inner product on $\Re$ as $\pair{\hat\Omega_1,\hat\Omega_2}=\frac{1}{2}\mathrm{tr}[\hat\Omega_1^T\hat\Omega_2]=\Omega_1\Omega_2$ for any $\Omega_1,\Omega_2\in\Re$. The matrix exponential is a local diffeomorphism from $\mathfrak{so}(2)$ to $\mathrm{SO}(2)$ given by
\begin{align*}
    \exp \hat\Omega = \begin{bmatrix} \cos\Omega & -\sin\Omega\\ \sin\Omega & \cos\Omega \end{bmatrix}.
\end{align*}
The kinematics equation for $R\in\mathrm{SO}(2)$ can be written in terms of a Lie algebra element as
\begin{align}
    \dot R = R\hat\Omega.\label{eqn:Rdotpp}
\end{align}

\textit{Continuous equations of motion:}
The Lagrangian for a planar pendulum $L:\mathrm{SO}(2)\times\mathfrak{so}(2)\rightarrow \Re$ can be written as
\begin{align}
    L(R,\hat\Omega) = \frac{1}{2}ml^2 \Omega^2 + mgl e_2^T R e_2= \frac{1}{2}ml^2 \pair{\hat\Omega,\hat\Omega} + mgl e_2^T R e_2,\label{eqn:Lpp}
\end{align}
where the constant $g\in\Re$ is the gravitational acceleration. The mass and the length of the pendulum are denoted by $m,l\in\Re$, respectively. The second expression is used to define a discrete Lagrangian later. We choose the bases of the inertial frame and the body-fixed frame such that the unit vector along the gravity direction in the inertial frame, and the unit vector along the pendulum axis in the body-fixed frame are represented by the same vector $e_2=[0;1]\in\Re^2$. Thus, for example, the hanging attitude is represented by $R=I_{2\times 2}$. Here, the rotation matrix $R\in\mathrm{SO}(2)$ represents the linear transformation from a representation of a vector in the body-fixed frame to the inertial frame.

Since the special orthogonal group is not a linear vector space, the expression for the variation should be carefully chosen. The infinitesimal variation of a rotation matrix $R\in\mathrm{SO}(2)$ can be written in terms of its Lie algebra element as
\begin{align}
    \delta R = \frac{d}{d\epsilon}\bigg|_{\epsilon=0} R \exp \epsilon\hat \eta = R\hat \eta,\label{eqn:delRpp}
\end{align}
where $\eta\in\Re$ so that $\hat\eta\in\mathfrak{so}(2)$.
The infinitesimal variation of the angular velocity is induced from this expression and \refeqn{Rdotpp} as
\begin{align}
    \delta \hat\Omega = \delta R^T \dot R + R^T \delta \dot R = -\hat\eta\hat\Omega + \hat\Omega\hat\eta + \hat{\dot\eta} = \hat{\dot\eta}.\label{eqn:delOmegapp}
\end{align}

Define the action integral to be $\mathfrak{G} = \int_{0}^T L(R,\hat\Omega)\,dt$. The infinitesimal variation of the action integral is obtained by using \refeqn{delRpp} and \refeqn{delOmegapp}. Hamilton's principle yields the following continuous equations of motion.
\begin{gather}
    \dot\Omega + \frac{g}{l} e_2^T R e_1=0,\\
    \dot R= R\hat\Omega.\label{eqn:dotR0}
\end{gather}
If we parameterize the rotation matrix as $R=\begin{bmatrix} \cos\theta & -\sin\theta \\ \sin\theta & \cos\theta\end{bmatrix}$, these equations are equivalent to
\begin{align}
    \ddot\theta + \frac{g}{l} \sin\theta =0.\label{eqn:ddottheta}
\end{align}

\textit{Discrete equations of motion:} We develop a Lie group variational integrator on $\mathrm{SO}(2)$. Similar to \refeqn{g1}, define $F_k\in\mathrm{SO}(2)$ such that
\begin{align}
    R_{k+1}=R_k F_k.\label{eqn:Rkppp}
\end{align}
Thus, $F_k = R_k^T R_{k+1}$ represents the relative update between two integration steps. If we find $F_k\in\mathrm{SO}(2)$, the orthogonal structure is preserved through \refeqn{Rkppp} since multiplication of orthogonal matrices is also orthogonal. This is a key idea of Lie group variational integrators.

Define the discrete Lagrangian $L_d:\mathrm{SO}(2)\times\mathrm{SO}(2)\rightarrow\Re$ to be
\begin{align}
    L_d(R_k,F_k)&=\frac{1}{2h}ml^2\pair{F_k-I_{2\times 2},F_k-I_{2\times 2}}+\frac{h}{2}mgle_2^T R_k e_2+\frac{h}{2}mgle_2^T R_{k+1} e_2,\nonumber\\
    & = \frac{1}{2h}ml^2\tr{I_{2\times 2}-F_k}+hmgle_2^T R_k e_2+\frac{h}{2}mgle_2^T R_{k+1} e_2,
\end{align}
which is obtained by an approximation $h\hat\Omega_k \simeq R_k^T (R_{k+1}-R_k)=F_k-I_{2\times 2}$, applied to the continuous Lagrangian given by \refeqn{Lpp}.

As for the continuous time case, expressions for the infinitesimal variations should be carefully chosen. The infinitesimal variation of a rotation matrix is the same as \refeqn{delRpp}, namely
\begin{align}
    \delta R_k = R_k \hat\eta_k,\label{eqn:delRkpp}
\end{align}
for $\eta_k\in\Re$, and the constrained variation of $F_k$ is obtained from \refeqn{Rkppp} as
\begin{align}
    \delta F_k = \delta R_k^T R_{k+1} + R_k^T \delta R_{k+1} = -\hat\eta_k F_k + F_k\hat\eta_{k+1}=F_k( \hat\eta_{k+1} - F_k^T \hat\eta_k F_k)=F_k(\hat\eta_{k+1}-\hat\eta_k),\label{eqn:delFkpp}
\end{align}
where we use the fact that $F\hat\eta F^T=\hat\eta$ for any $F\in\mathrm{SO}(2)$ and $\hat\eta\in\mathfrak{so}(2)$.

Define an action sum $\mathfrak{G}_d:\mathrm{SO}(2)^{N+1}\rightarrow\Re$ as $\mathfrak{G}_d = \sum_{k=0}^{N-1} L_d(R_k,F_k)$. Using \refeqn{delRkpp} and \refeqn{delFkpp}, the variation of the action sum is written as
\begin{align*}
    \delta\mathfrak{G}_d = \sum_{k=1}^{N-1}\pairlr{\frac{1}{2h}ml^2(F_{k-1}-F_{k-1}^T)-\frac{1}{2h}(F_{k}-F_{k}^T)-hmgl\widehat{e_2^T R_k e_1},\;\hat\eta_{k}}.
\end{align*}
From discrete Hamilton's principle, $\delta\mathfrak{G}_d=0$ for any $\hat\eta_k$. Thus, we obtain the Lie group variational integrator on $\mathrm{SO}(2)$ as
\begin{gather}
    (F_{k}-F_{k}^T)-(F_{k+1}-F_{k+1}^T)-\frac{2 h^2 g}{l}\widehat{e_2^T R_{k+1} e_1}=0,\label{eqn:findf0}\\
    R_{k+1}=R_k F_k.\label{eqn:Rkp0}
\end{gather}
For a given $(R_k,F_k)$ and $R_{k+1}=R_kF_k$, \refeqn{findf0} is solved to find $F_{k+1}$. This yields a discrete map $(R_k,F_k)\mapsto(R_{k+1},F_{k+1})$. If we parameterize the rotation matrices $R$ and $F$ with $\theta$ and $\Delta\theta$ and if we assume that $\Delta\theta \ll 1$, these equations are equivalent to
\begin{align*}
    \frac{1}{h}(\theta_{k+1}-2\theta_k +\theta_k) + \frac{hg}{l} \sin \theta_{k}=0.
\end{align*}

The discrete version of the Legendre transformation provides the discrete Hamiltonian map as follows.
\begin{gather}
    F_k - F_k^T=2h\hat\Omega -\frac{h^2 g}{l}\widehat{e_2^T R_k e_1},\label{eqn:findfpp}\\
    R_{k+1}=R_k F_k,\label{eqn:Rkp1pp}\\
    \Omega_{k+1} = \Omega_k -\frac{hg}{2l} e_2^T R_k e_1-\frac{hg}{2l}e_2^T R_{k+1} e_1.\label{eqn:mukppp}
\end{gather}
For a given $(R_k,\Omega_k)$, we solve \refeqn{findfpp} to obtain $F_k$. Using this, $(R_{k+1},\Omega_{k+1})$ is obtained from \refeqn{Rkp1pp} and \refeqn{mukppp}. This yields a discrete map $(R_k,\Omega_k)\mapsto(R_{k+1},\Omega_{k+1})$.

\textit{Numerical example:} We compare the computational properties of the discrete equations of motion given by \refeqn{findfpp}--\refeqn{mukppp} with a 4(5)-th order variable step size Runge-Kutta method. We choose $m=1\,\mathrm{kg}$, $l=9.81\,\mathrm{m}$. The initial conditions are $\theta_0=\pi/2\,\mathrm{rad}$, $\Omega=0$, and the total energy is $E=0\,\mathrm{Nm}$. The simulation time is $1000\,\mathrm{sec}$, and the step-size $h=0.03$ of the discrete equations of motion is chosen such that the CPU times are identical. \reffig{Epp} shows the computed total energy for both methods. The variational integrator preserves the total energy well. There is no drift in the computed total energy, and the mean variation is $1.0835\times 10^{-2}\,\mathrm{Nm}$. But, there is a notable dissipation of the computed total energy for the Runge-Kutta method. Note that the computed total energy would further decrease as the simulation time increases.
\end{example}

\begin{figure}
\subfigure[Trajectory of pendulum]%
    {\includegraphics[width=0.24\textwidth]{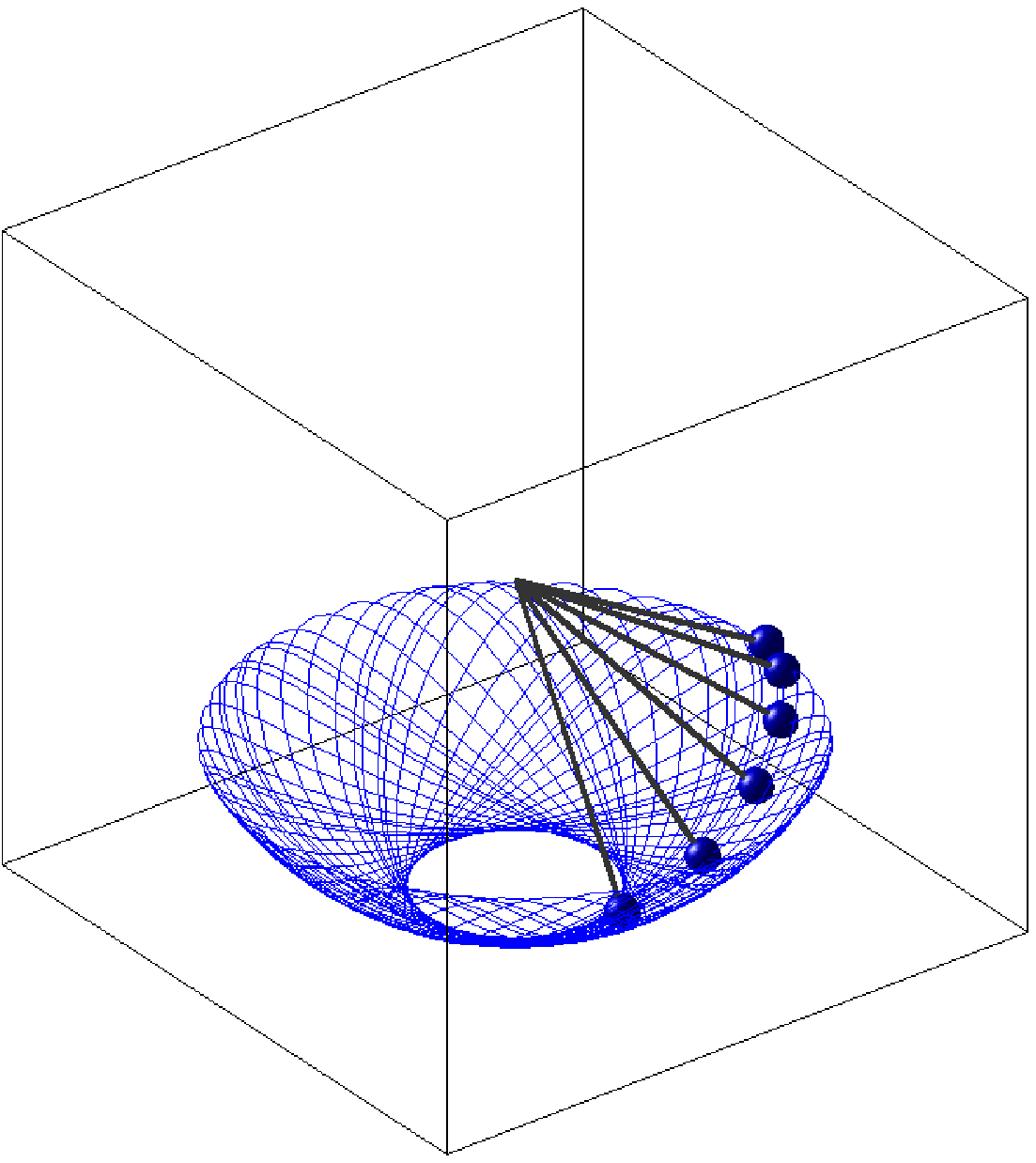}}
\hspace{0.002\textwidth}
\subfigure[Total energy]%
    {\includegraphics[width=0.35\textwidth]{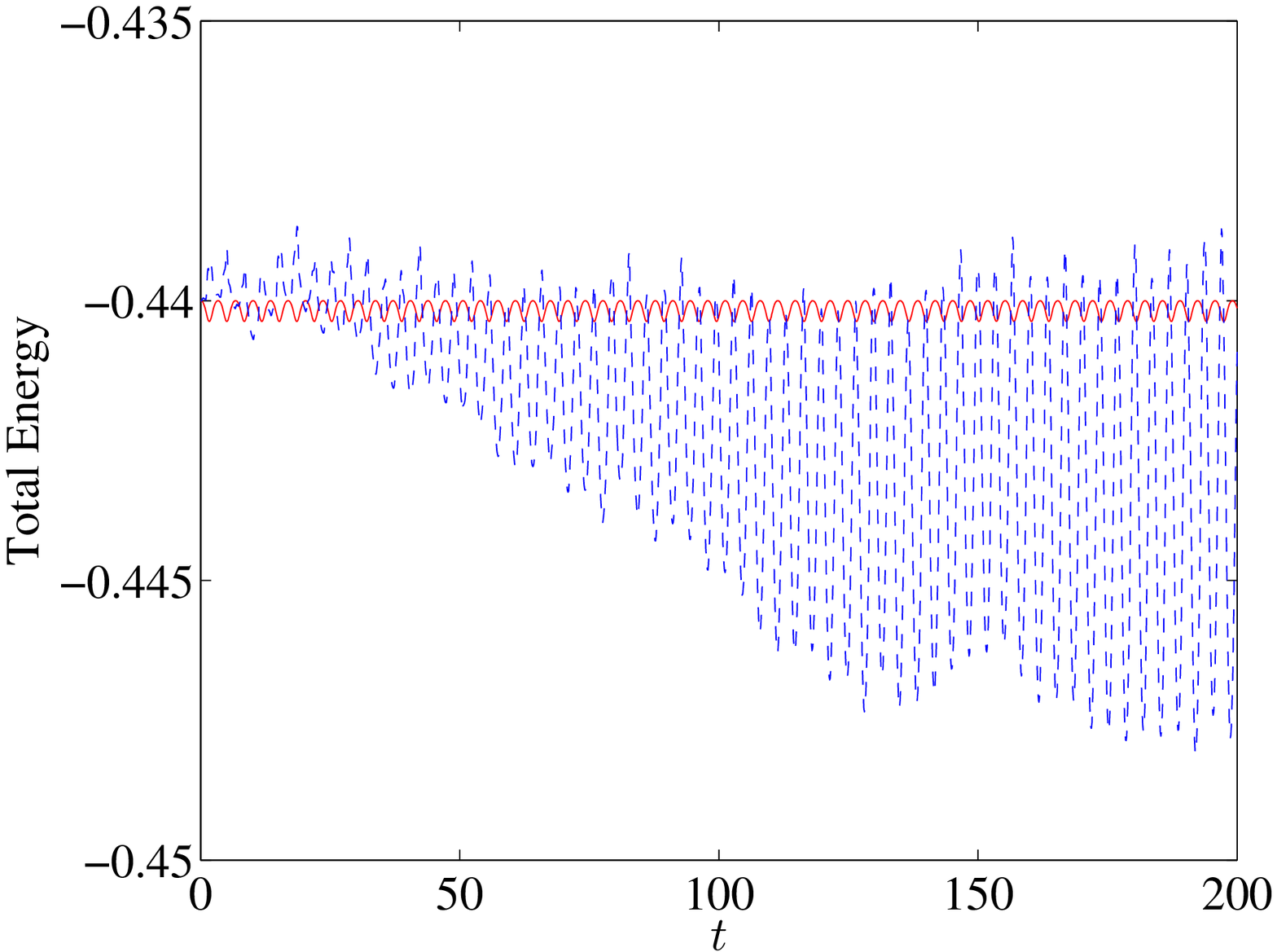}}
\subfigure[Unit length error $| q_k\cdot q_k-1|$]%
    {\includegraphics[width=0.35\textwidth]{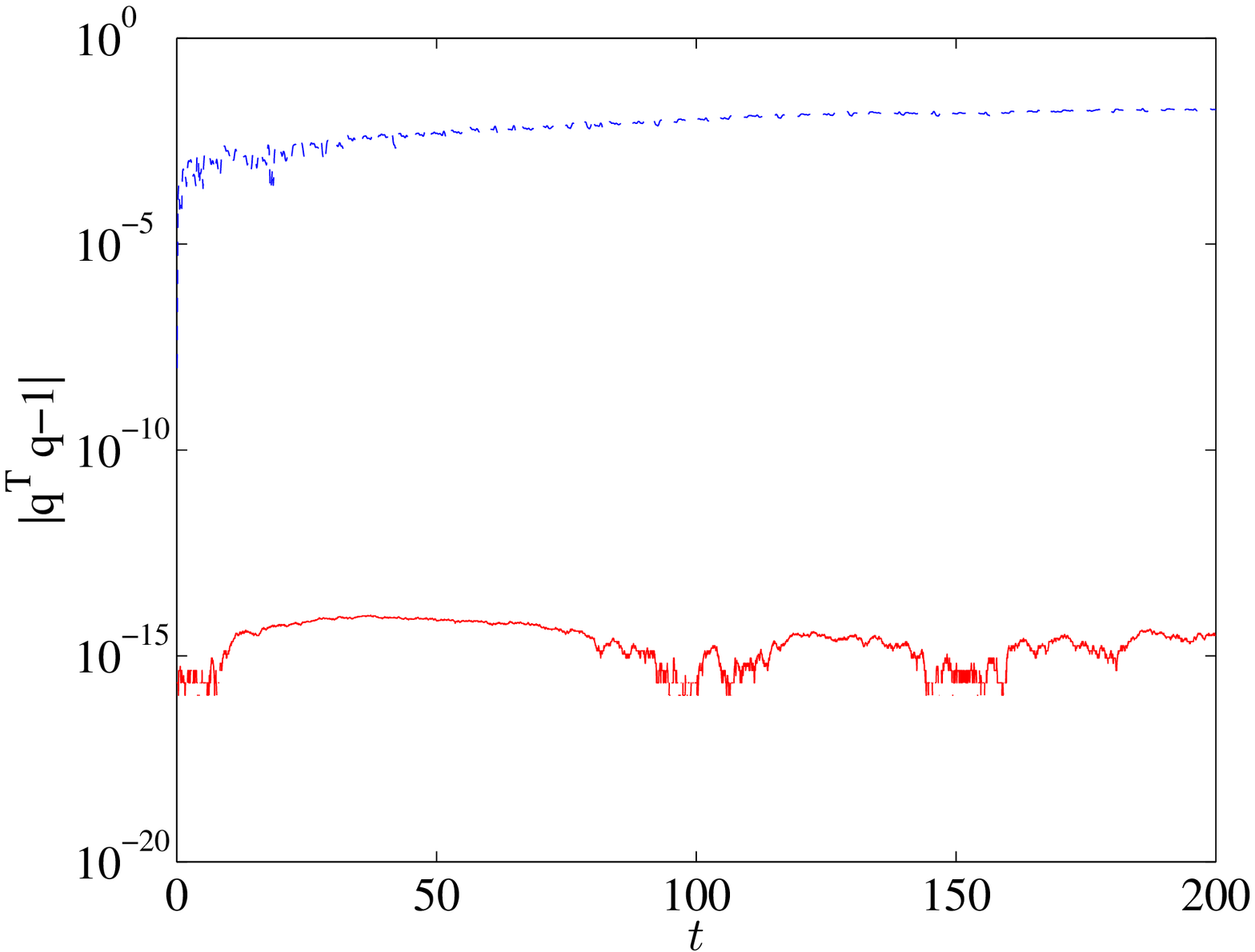}}
\caption{Numerical simulation of a spherical pendulum (RK45: blue, dotted, VI: red, solid)}\label{fig:sp}
\end{figure}

\begin{example}[\textbf{Spherical Pendulum}]\label{ex:sp}
A spherical pendulum is a mass particle connected to a frictionless, two degree-of-freedom pivot by a rigid massless link. The mass particle acts under a uniform gravitational potential. The configuration space is the two-sphere $\S^2=\braces{q\in\Re^3\,|\,\norm{q}=1}$. It is common to parameterize the two-sphere by two angles, but this description of the spherical pendulum has a singularity. Any trajectory near to singularity causes numerical ill-conditioning. Furthermore, this leads to complicated expressions involving trigonometric functions.

Here we develop equations of motion for a spherical pendulum using the global structure of the two-sphere without parameterization. In the previous example, we develop equations of motion for a planar pendulum using the fact that the one-sphere $\S^1$ is diffeomorphic to the special orthogonal group $\mathrm{SO}(2)$. But, the two-sphere is not diffeomorphic to a Lie group. Instead, there exists a natural Lie group action on the two-sphere. That is the 3-dimensional special orthogonal group $\SO$, a group of $3\times 3$ orthogonal matrices with determinant 1, i.e. $\SO=\braces{R\in\Re^{3\times 3}\,|\,R^TR=I_{3\times 3},\mathrm{det}[R]=1}$. The special orthogonal group $\SO$ acts on the two-sphere in a transitive way; for any $q_1,q_2\in\S^2$, there exists a $R\in\SO$ such that $q_2=Rq_1$. Thus, the discrete update for the two-sphere can be expressed in terms of a rotation matrix as \refeqn{Rkppp}. This is a key idea to develop a discrete equations of  motion for a spherical pendulum.

\textit{Continuous equations of motion:}
Let $q\in\S^2$ be a unit vector from the pivot point to the point mass. The Lagrangian for a spherical pendulum can be written as
\begin{align}
    L(q,\dot q) = \frac{1}{2} ml^2 \dot q \cdot \dot q + mgl e_3 \cdot q,
\end{align}
where the gravity direction is assumed to be $e_3=[0;0;1]\in\Re^3$. The mass and the length of the pendulum are denoted by $m,l\in\Re$, respectively.
The infinitesimal variation of the unit vector $q$ can be written in terms of the vector cross product as
\begin{align}
    \delta q = \xi \times q\label{eqn:delqsp},
\end{align}
where $\xi\in\Re^3$ is constrained to be orthogonal to the unit vector, i.e. $\xi\cdot q=0$. Using this expression for the infinitesimal variation, Hamilton's principle yields the following continuous equations of motion.
\begin{align}
    \ddot q + (\dot q\cdot \dot q)q  + \frac{g}{l} (q\times (q\times e_3))=0.
\end{align}
Since $\dot q= \omega\times q$ for some angular velocity $\omega\in\Re^3$ satisfying $\omega\cdot q=0$, this can also be equivalently written as
\begin{gather}
    \dot\omega  - \frac{g}{l} q\times e_3=0,\label{eqn:omegadotsp}\\
    \dot q = \omega\times q.\label{eqn:qdotsp}
\end{gather}
These are global equations of motion for a spherical pendulum; these are much simpler than the equations expressed in term of angles, and they have no singularity.

\textit{Discrete equations of motion:}
We develop a variational integrator for the spherical pendulum defined on $\S^2$. Since the special orthogonal group acts on the two-sphere transitively, we can define the discrete update map for the unit vector as
\begin{align}
    q_{k+1} = F_k q_k\label{eqn:qkp0sp}
\end{align}
for $F_k\in\SO$. Define a discrete Lagrangian $L_d:\S^2\times\S^2\rightarrow\Re$ to be
\begin{align*}
    L_d (q_k,q_{k+1}) & = \frac{1}{2h}ml^2 (q_{k+1}-q_k)\cdot (q_{k+1}-q_k)+\frac{h}{2}mgl e_3\cdot q_k+\frac{h}{2}mgl e_3\cdot q_{k+1}.
\end{align*}
The variation of $q_k$ is the same as \refeqn{delqsp}, namely
\begin{align}
    \delta q_k = \xi_k \times q_k
\end{align}
for $\xi_k\in\Re^3$ with a constraint $\xi_k\cdot q_k=0$. Using this discrete Lagrangian and the expression for the variation, discrete Hamilton's principle yields the following discrete equations of motion for a spherical pendulum.
\begin{gather}
    q_{k+1}=\parenth{h\omega_k + \frac{h^2g}{2l} q_k \times e_3}\times q_k + \parenth{1-\norm{h\omega_k + \frac{h^2g}{2l} q_k \times e_3}^2}^{1/2}q_k,\label{eqn:qkpesp}\\
    \omega_{k+1} = \omega_k+\frac{hg}{2l} q_{k} \times e_3  +\frac{hg}{2l} q_{k+1} \times e_3,\label{eqn:omegakpesp}
\end{gather}
Since an explicit solution for $F_k\in\SO$ can be obtained in this case, the rotation matrix $F_k$ does not appear in the equations of motion. This variational integrator on $\S^2$ exactly preserves the unit length of $q_k$, the constraint $q_k\cdot \omega_k=0$, and the third component of the angular velocity $\omega_k\cdot e_3$ which is conserved since gravity exerts no moment along the gravity direction $e_3$.

\textit{Numerical example:}
We compare the computational properties of the discrete equations of motion given by \refeqn{qkpesp} and \refeqn{omegakpesp} with a 4(5)-th order variable step size Runge-Kutta method for \refeqn{omegadotsp} and \refeqn{qdotsp}. We choose $m=1\,\mathrm{kg}$, $l=9.81\,\mathrm{m}$. The initial conditions are $q_0=[\frac{\sqrt{3}}{2},0,\frac{1}{2}]$, $\omega_0=0.1[\sqrt{3},0,3]\,\mathrm{rad/sec}$, and the total energy is $E=-0.44\,\mathrm{Nm}$. The simulation time is $200\,\mathrm{sec}$, and the step-size $h=0.05$ of the discrete equations of motion is chosen such that the CPU times are identical. \reffig{sp} shows the computed total energy and the unit length errors. The variational integrator preserves the total energy and the structure of the two-sphere well. The mean total energy deviation is $1.5460\times 10^{-4}\,\mathrm{Nm}$, and the mean unit length error is $3.2476\times 10^{-15}$. But, there is a notable dissipation of the computed total energy for the Runge-Kutta method. The Runge-Kutta method also fails to preserve the structure of the two-sphere. The mean unit length error is $1.0164\times 10^{-2}$.
\end{example}

\begin{example}[\textbf{Rigid Body in a Potential Field}]\label{ex:rb}
Consider a rigid body under a potential field that is dependent on the position and the attitude of the body. The configuration space is the special Euclidean group, which is a semi-direct product of the special orthogonal group and Euclidean space, i.e. $\SE=\SO\,\textcircled{s}\,\Re^3$.

\textit{Continuous equations of motion:}
The equations of motion for a rigid body can be developed either from Hamilton's principle (see \citepri{CMA07}) in a similar way as Example \ref{ex:pp}, or directly from the generalized discrete Euler--Poincar\'{e} equation given at \refeqn{DEP}. Here, we summarize the results. Let $m\in\Re$ and $J\in\Re^{3\times 3}$ be the mass and the moment of inertia matrix of a rigid body. For $(R,x)\in\SE$, the linear transformation from the body-fixed frame to the inertial frame is denoted by the rotation matrix $R\in\SO$, and the position of the mass center in the inertial frame is denoted by a vector $x\in\Re^3$. The vectors $\Omega,v\in\Re^3$ are the angular velocity in the body-fixed frame, and the translational velocity in the inertial frame, respectively. Suppose that the rigid body acts under a configuration-dependent potential $U:\SE\rightarrow\Re$. The continuous equations of motion for the rigid body can be written as
\begin{gather}
\dot{R}=R\hat\Omega,\label{eqn:Rdot}\\
\dot{x}=v,\\
J\dot{\Omega}+\Omega\times J\Omega = M,\\
m\dot{v}=-\deriv{U}{x},\label{eqn:vdot}
\end{gather}
where the hat map $\hat\cdot$ is an isomorphism from $\Re^3$ to $3\times 3$ skew-symmetric matrices $\so$, defined such that $\hat x y =x\times y$ for any $x,y\in\Re^3$. The moment due to the potential $M\in\Re^3$ is obtained by the following relationship.
\begin{align}
\hat M & = \deriv{U}{R}^TR-R^T\deriv{U}{R}.
\end{align}
The matrix $\deriv{U}{R}\in\Re^{3\times 3}$ is defined such that $[\deriv{U}{R}]_{ij}=\deriv{U}{[R]_{ij}}$ for $i,j\in\braces{1,2,3}$, where the $i,j$-th element of a matrix is denoted by $[\cdot]_{ij}$.

\textit{Discrete equations of motion:}
The corresponding discrete equations of motion are given by
\begin{gather}
h \widehat{J\Omega_{k}}+\frac{h^2}{2}\hat M_{k}=F_{k}J_{d}-J_{d}F_{k}^T,\label{eqn:findF}\\
R_{k+1}=R_{k}F_{k},\label{eqn:Rkp}\\
x_{k+1} = x_{k} + hv_k -\frac{h^2}{2m}\deriv{U_{{k}}}{x_{k}},\label{eqn:xkp}\\
J\Omega_{{k+1}}=F_{k}^T J\Omega_{k}+\frac{h}{2}F_{k}^TM_{k}+\frac{h}{2}M_{{k+1}}
\label{eqn:Omegakp},\\
m v_{k+1}=m v_k-\frac{h}{2}\deriv{U_k}{x_k}
-\frac{h}{2}\deriv{U_{k+1}}{x_{k+1}},\label{eqn:vkp}
\end{gather}
where $J_d\in\Re^{3\times 3}$ is a non-standard moment of inertia matrix defined as $J_d=\frac{1}{2}\mathrm{tr}[J]I_{3\times 3}-J$. For a given $(R_k,x_k,\Omega_k,v_k)$, we solve the implicit equation \refeqn{findF} to find $F_k\in\SO$. Then, the configuration at the next step $R_{k+1},x_{k+1}$ is obtained by \refeqn{Rkp} and \refeqn{xkp}, and the moment and force $M_{k+1},\deriv{U_{k+1}}{x_{k+1}}$ can be computed. Velocities $\Omega_{k+1},v_{k+1}$  are obtained from \refeqn{Omegakp} and \refeqn{vkp}. This
defines a discrete flow map, $(R_k,x_k,\Omega_k,v_k)\mapsto(R_{k+1},x_{k+1},\Omega_{k+1},v_{k+1})$,
and this process can be repeated. This Lie group variational integrator on $\SE$ can be generalized to multiple rigid bodies acting under their mutual gravitational potential (see \citepri{CMA07}).

\begin{figure}
\subfigure[Trajectory of rigid body]%
    {\includegraphics[width=0.30\textwidth]{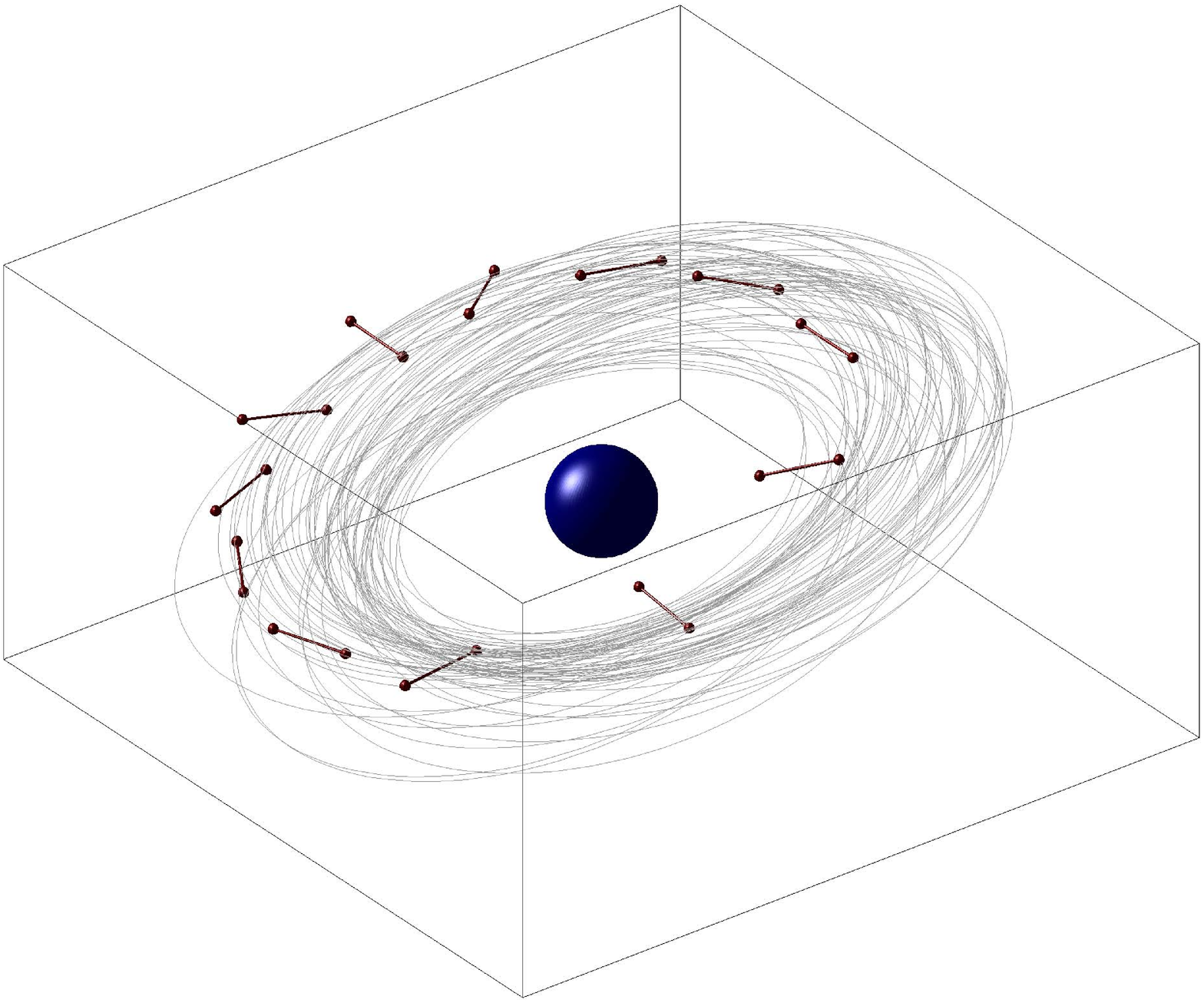}\label{fig:trajrb}}
\hspace{0.002\textwidth}
\subfigure[Total energy]%
    {\includegraphics[width=0.33\textwidth]{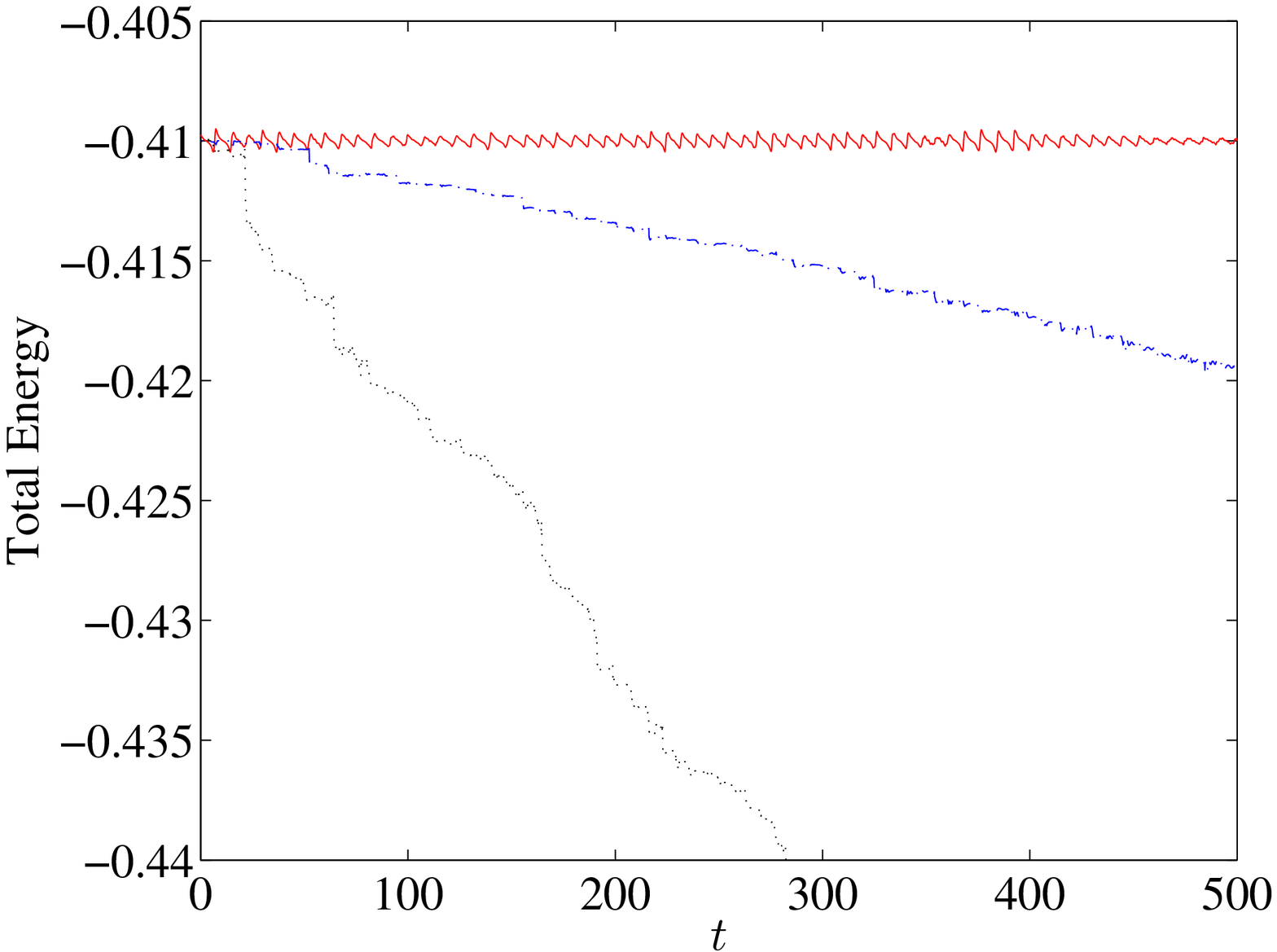}\label{fig:Erb}}
\subfigure[Rotation matrix error $\norm{I-R^T R}$]%
    {\includegraphics[width=0.33\textwidth]{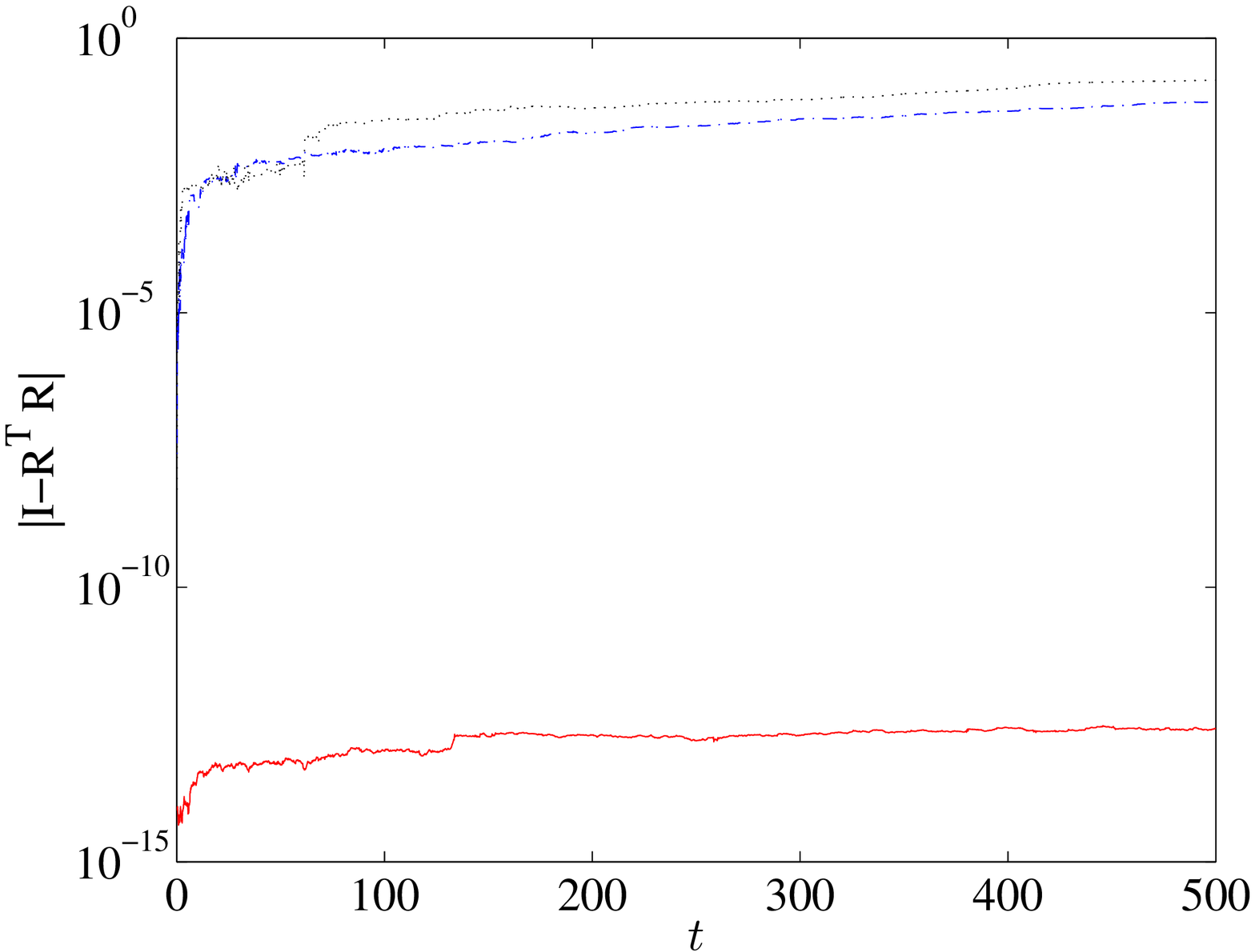}\label{fig:errRrb}}
\caption{Numerical simulation of a dumbbell rigid body (LGVI: red, solid, RK45 with rotation matrices: blue, dash-dotted, RK45 with quaternions: black, dotted)}\label{fig:rb}
\end{figure}

\textit{Numerical example:}
We compare the computational properties of the discrete equations of motion given by \refeqn{findF}--\refeqn{vkp} with a 4(5)-th order variable step size Runge-Kutta method for \refeqn{Rdot}--\refeqn{vdot}. In addition, we compute the attitude dynamics using quaternions on the unit three-sphere $\S^3$. The attitude kinematics equation \refeqn{Rdot} is rewritten in terms of quaternions, and the corresponding equations are integrated by the same Runge-Kutta method.

We choose a dumbbell spacecraft, that is two spheres connected by a rigid massless rod, acting under a central gravitational potential. The resulting system is a restricted full two body problem. The dumbbell spacecraft model has an analytic expression for the gravitational potential, resulting in a nontrivial coupling between the attitude dynamics and the orbital dynamics.

As shown in \reffig{trajrb}, the initial conditions are chosen such that the resulting motion is a near-circular orbit combined with a rotational motion. \reffig{Erb} and \reffig{errRrb} show the computed total energy and the orthogonality errors of the rotation matrix. The Lie group variational integrator preserves the total energy and the Lie group structure of $\SO$. The mean total energy deviation is $2.5983\times 10^{-4}$, and the mean orthogonality error is $1.8553\times 10^{-13}$. But, there is a notable dissipation of the computed total energy and the orthogonality error for the Runge-Kutta method. The mean orthogonality errors for the Runge-Kutta method are $0.0287$ and $0.0753$, respectively, using kinematics equation with rotation matrices, and using the kinematics equation with quaternions. Thus, the attitude of the rigid body is not accurately computed for Runge-Kutta methods. It is interesting to see that the Runge-Kutta method with quaternions, which is generally assumed to have better computational properties than the kinematics equation with rotation matrices, has larger total energy error and orthogonality error. Since the unit length of the quaternion vector is not preserved in the numerical computations, orthogonality errors arise when converted to a rotation matrix. This suggests that it is critical to preserve the structure of $\SO$ in order to study the global characteristics of the rigid body dynamics.

The dynamics of a rigid body is characterized by a Hamiltonian system on a Lie group. The Lie group variational integrator has a desirable property that both the symplectic structure and the Lie group structure of the rigid body dynamics are preserved concurrently. More explicitly, the computational properties of the Lie group variational integrator is compared with a symplectic integrator that does not preserve the Lie group structure, and a Lie group method that does not preserve the symplectic structure (see \citepri{CMDA07}). It is shown that the Lie group variational integrator has a substantial superiority in terms of numerical accuracy and efficiency. Due to these computational advantages, the Lie group variational integrator has been used to study the dynamics of the binary near-Earth asteroid 66391 (1999 $KW_4$) in joint work between the University of Michigan and the Jet Propulsion Laboratory, NASA (see \citepri{Sch.Sci06}).

\end{example}

\section{Optimal Control of Discrete Lagrangian and Hamiltonian System}
Optimal control problems involve finding a control input such that a certain optimality objective is achieved under prescribed constraints. An optimal control problem that minimizes a performance index is described by a set of differential equations, which can be derived using Pontryagin's maximum principle. The equations of motion for a system are constrained by Lagrange multipliers, and necessary conditions for optimality is obtained by the calculus of variation. The solution for the corresponding two point boundary value problem provides the optimal control input. Alternatively, a sub-optimal control law is obtained by approximating the control input history with finite data points.

Discrete optimal control problems involve finding a control input for a given system described by discrete Lagrangian and Hamiltonian mechanics. The control inputs are parameterized by their values at each discrete time step, and the discrete equations of motion are derived from the discrete Lagrange-d'Alembert principle (\citepri{KaMaOrWe2000}),
\[\delta \sum_{k=0}^{N-1} L_d(q_k,q_{k+1})=\sum_{k=0}^{N-1} [F_d^{-}(q_k,q_{k+1})\cdot\delta q_k+F_d^+(q_k,q_{k+1})\cdot\delta q_{k+1}],\]
which modifies the discrete Hamilton's principle by taking into account the virtual work of the external forces.
Discrete optimal control is in contrast to traditional techniques such as collocation, wherein the continuous equations of motion are imposed as constraints at a set of collocation points, since this approach induces constraints on the configuration at each discrete timestep.

Any optimal control algorithm can be applied to discrete Lagrangian or Hamiltonian system. For an indirect method, our approach to a discrete optimal control problem can be considered as a multiple stage variational problem. The discrete equations of motion are derived by the discrete variational principle. The corresponding variational integrators are imposed as dynamic constraints to be satisfied by using Lagrange multipliers, and necessary conditions for optimality, expressed as  discrete equations on multipliers, are obtained from a variational principle. For a direct method, control inputs can be optimized by using parameter optimization tools such as a sequential quadratic programming. The discrete optimal control can be characterized by discretizing the optimal control problem from the problem formulation stage.

This method has substantial computational advantage to find an optimal control law. As discussed in the previous section, the discrete dynamics are more faithful to the continuous
equations of motion, and consequently more accurate solutions to the optimal control problem are obtained. The external control inputs break the Lagrangian and Hamiltonian system structure. For example, the total energy is not conserved for a controlled mechanical system. But, the computational superiority of the discrete mechanics still holds for controlled systems. It has been shown that the discrete dynamics is more reliable even for controlled system as it computes the energy dissipation rate of controlled systems more accurately (see \citepri{MaWe2001}). For example, this feature is extremely important in computing accurate optimal trajectories for long term spacecraft attitude maneuvers using low energy control inputs.

The discrete dynamics does not only provide an accurate optimal control input, but also enables us to find it efficiently. For the indirect optimal control approach, optimal solutions are usually sensitive to a small variation of multipliers. This causes difficulties, such as numerical ill-conditioning, to solve the necessary conditions for optimality expressed as a two point boundary value problem. Sensitivity derivatives along the discrete necessary conditions do not have numerical dissipation caused by conventional numerical integration schemes. Thus, they are numerically more robust, and the necessary conditions can be solved computationally efficiently. For the direct optimal control approach, optimal control inputs can be obtained by using larger discrete step size, which requires less computational load.

We illustrate the basis properties of the discrete optimal control using optimal control problems for the spherical pendulum and the rigid body model presented in the previous section.

\begin{example}[\textbf{Optimal Control of a Spherical Pendulum}]
We study an optimal control problem for the a spherical pendulum described in Example \ref{ex:sp}.
We assume that an external control moment $u\in\Re^3$ acts on the pendulum. Control inputs are parameterized by their values at each time step, and the discrete equations of motion are modified to include the effects of the external control inputs by using discrete Lagrange-d'Alembert principle. Since the component of the control moment that is parallel to the direction along the pendulum has no effect, we parameterize the control input as $u_k=q_k\times w_k$ for $w_k\in\Re^3$.

The objective is to transfer the pendulum from a given initial configuration
$(q_0,\omega_0)$ to a desired configuration
$(q^d,\omega^d)$ during a fixed maneuver time
$N$, while minimizing the square of the weighted $l_2$ norm
of the control moments.
\begin{gather*}
\min_{w_{k}} \mathcal{J}=\sum_{k=0}^{N}
\frac{h}{2}u_k^T u_k = \sum_{k=0}^{N}\frac{h}{2} (q_k\times w_k)^T (q_k\times w_k).
\end{gather*}

We solve this optimal control problem by using a direct numerical optimization method. The terminal boundary condition is imposed as an equality constraint, and the $3N$ control input parameters $\braces{w_k}_{k=0}^N$ are numerically optimized using sequential quadratic programming. This method is referred to as a DMOC (Discrete Mechanics and Optimal Control) approach (see \citepri{JuMaOb2005}).

\begin{figure}
\subfigure[Pendulum position $q$]%
    {\includegraphics[width=0.302\textwidth]{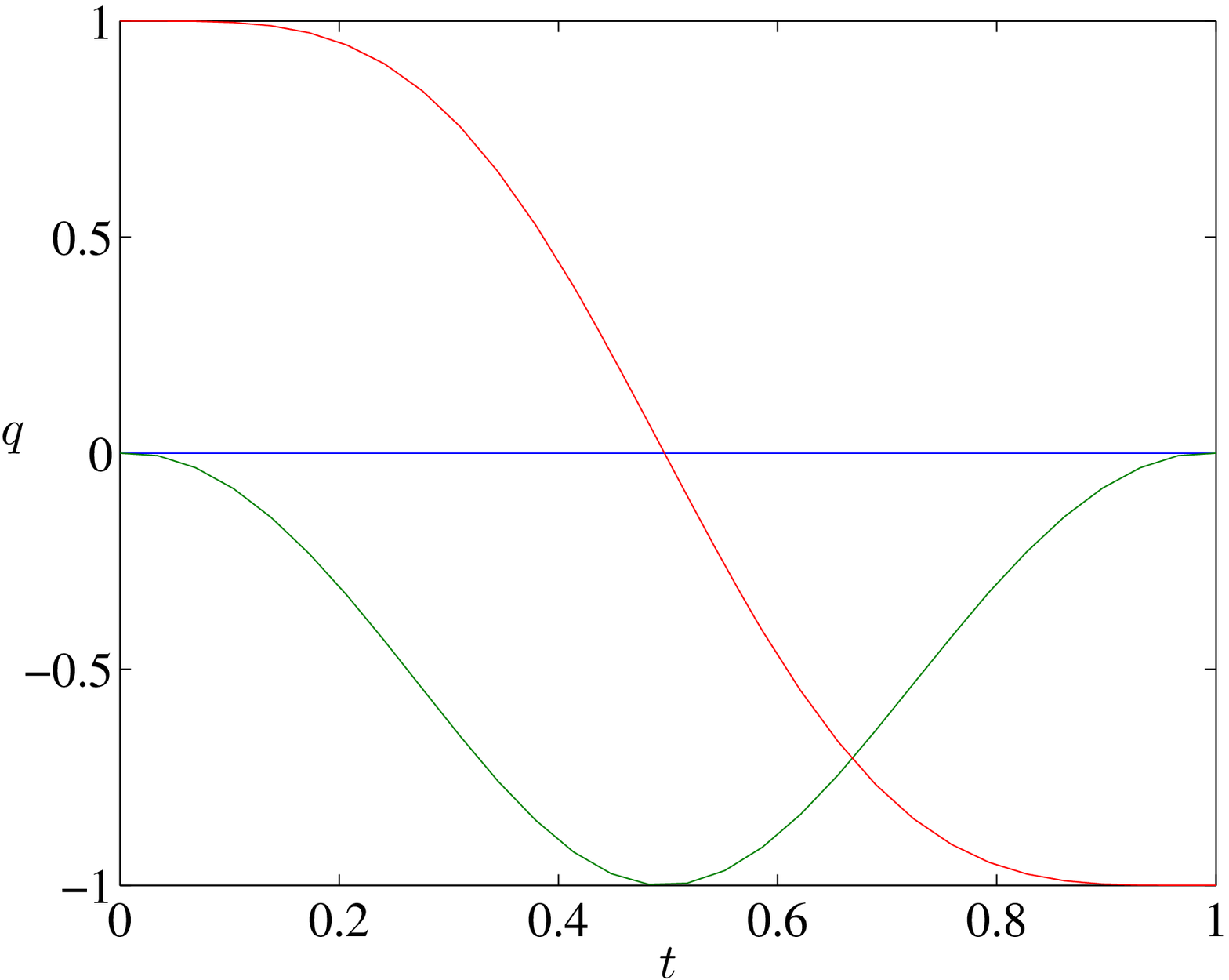}\label{fig:optspq}}
\hspace{0.002\textwidth}
\subfigure[Angular velocity $\omega$]%
    {\includegraphics[width=0.295\textwidth]{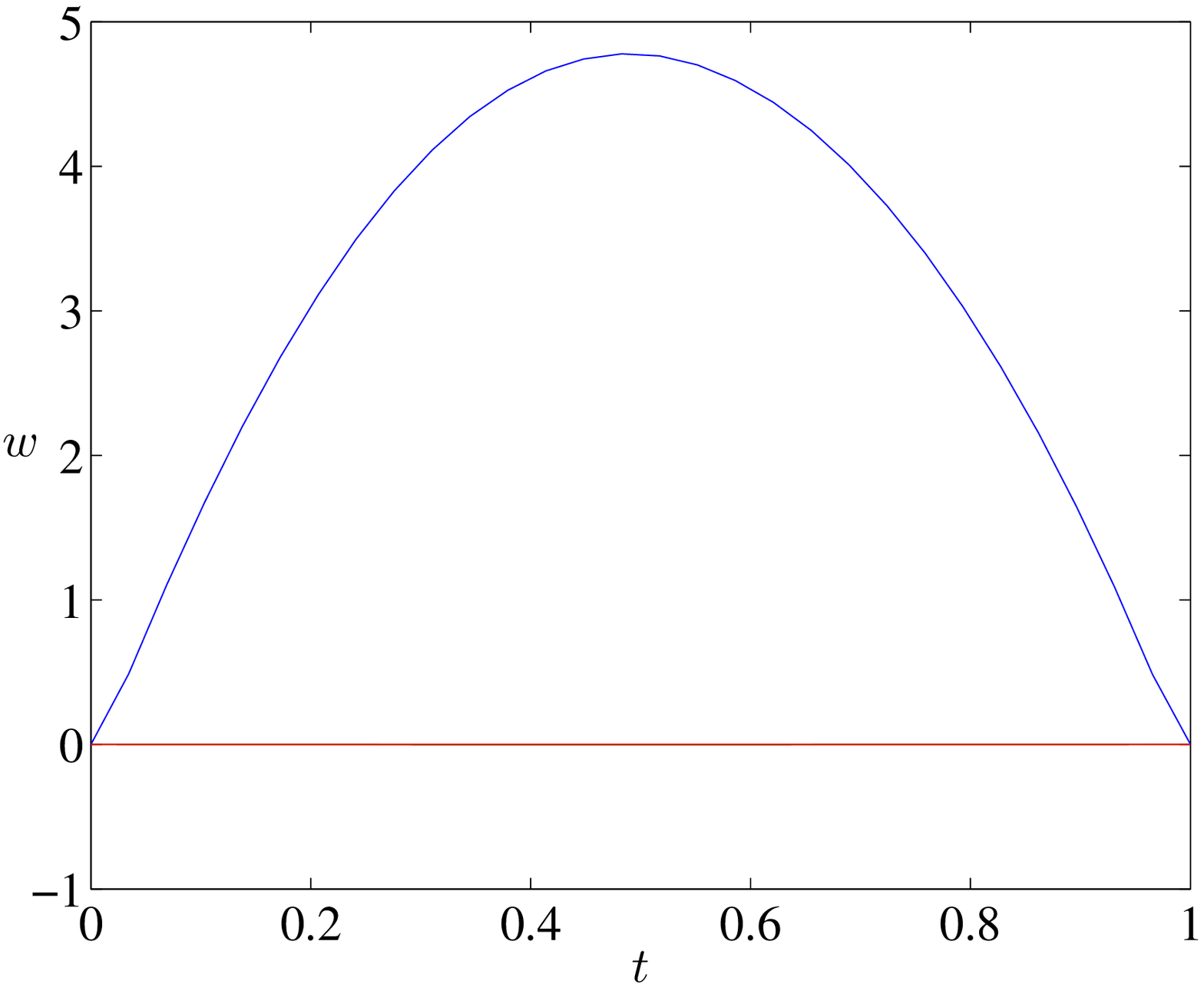}\label{fig:optspw}}
\hspace{0.002\textwidth}
\subfigure[Control moment $u$]%
    {\includegraphics[width=0.30\textwidth]{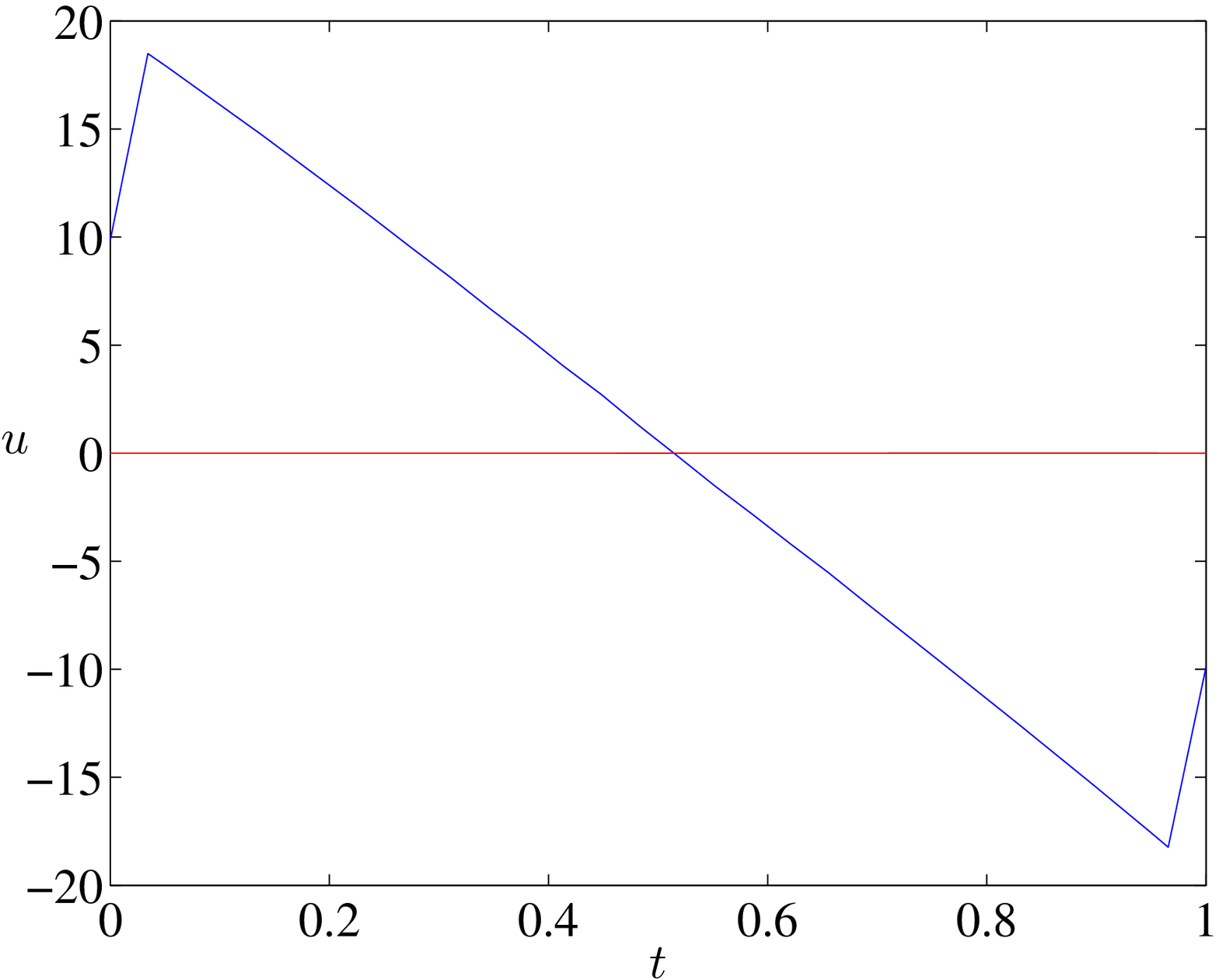}\label{fig:optspu}}
\caption{Optimal control of a spherical pendulum}\label{fig:optsp}
\end{figure}

\reffig{optsp} shows a optimal solution transferring the spherical pendulum from a hanging configuration given by $(q_0,\omega_0)=(e_3,0_{3\times 1})$ to an inverted configuration $(q^d,\omega^d)=(-e_3,0_{3\times 1})$ during 1 second. The time step size is $h=0.033$. Experiment has shown that the DMOC approach can compute optimal solutions using larger step size than typical Runge-Kutta methods, and consequently, it requires less computational load. In this case, using a second-order accurate Runge-Kutta method, the same optimization code fails while giving error messages of inaccurate and singular gradient computations. It is presumed that the unit length errors of the Runge-Kutta method, shown in Example \ref{ex:sp}, cause numerical instabilities for a finite-difference gradient computations required for the sequential quadratic programming algorithm.

\end{example}

\begin{example}[\textbf{Optimal Control of a Rigid Body in a Potential Field}]
We study an optimal control problem of a rigid body using a dumbbell spacecraft model described in Example \ref{ex:rb} (see \citepri{CDC06.opt} for detail). We assume that external control forces $u^f\in\Re^3$, and control moment $u^m\in\Re^3$ act on the dumbbell spacecraft. Control inputs are parameterized by their values at each time step, and the Lie group variational integrators are modified to include the effects of the external control inputs by using discrete Lagrange-d'Alembert principle.

The objective is to transfer the dumbbell from a given initial configuration
$(R_0,x_0,\Omega_0,v_0)$ to a desired configuration
$(R^d,x^d,\Omega^d,v^d)$ during a fixed maneuver time
$N$, while minimizing the square of the $l_2$ norm
of the control inputs.
\begin{gather*}
\min_{u_{k+1}} \mathcal{J}=\sum_{k=0}^{N-1}
\frac{h}{2}(u^f_{k+1})^TW_fu^f_{k+1}+
\frac{h}{2}(u^m_{k+1})^TW_mu^m_{k+1},
\end{gather*}
where $W_f,W_m\in\Re^{3\times 3}$ are symmetric positive definite
matrices.  Here we use a modified version of the discrete equations
of motion with first order accuracy, as it yields a compact
form for the necessary conditions.

\textit{Necessary conditions for optimality:}
We solve this optimal control problem by using an indirect optimization method, where necessary conditions for optimality are derived using variational arguments, and a solution of the corresponding two-point boundary value problem provides the optimal control. This approach is common in the optimal control literatures; here the optimal control problem is discretized at the problem formulation level using the Lie group variational integrator presented in Section \ref{sec:dm}.
\begin{align*}
\mathcal{J}_a =
\sum_{k=0}^{N-1}&\frac{h}{2}(u^f_{k+1})^TW^fu^f_{k+1}+
\frac{h}{2}(u^m_{k+1})^TW^mu^m_{k+1}\nonumber\\
& +\lambda_k^{1,T}\braces{-x_{k+1}+x_k+h v_k} +\lambda_k^{2,T}\braces{-m v_{k+1} + mv_k-h\deriv{U_{k+1}}{x_{k+1}}+hu^f_{k+1}}\nonumber\\
& +\lambda_k^{3,T}\parenth{\mathrm{logm}(F_k-R_{k}^TR_{k+1})}^{\vee} +\lambda_k^{4,T}\braces{-J\Omega_{k+1} + F_k^T J\Omega_k +
h\parenth{M_{k+1}+u_{k+1}^m}},
\end{align*}
where $\lambda_k^1,\lambda_k^2,\lambda_k^3,\lambda_k^4\in \Re^3$ are Lagrange multipliers. The matrix logarithm is denoted by $\mathrm{logm}:\SO\rightarrow\so$ and the vee map $\vee:\so\rightarrow\Re^3$ is the inverse of the hat map introduced in Example \ref{ex:rb}. The logarithm form of \refeqn{Rkp} is used, and the
constraint \refeqn{findF} is considered implicitly using constrained variations. Using  similar expressions for the variation of the rotation matrix and the angular velocity given in \refeqn{delRpp} and \refeqn{delOmegapp}, the infinitesimal variation can be written as
\begin{align*}
\delta\mathcal{J}_a & = \sum_{k=1}^{N-1} h\delta
u_{k}^{f,T}\braces{W_fu^f_{k}+\lambda_{k-1}^2}
+h\delta u_{k}^{m,T}\braces{W_mu^m_{k}+\lambda_{k-1}^4}
+z_k^T\braces{-\lambda_{k-1}+A_k^T \lambda_k},
\end{align*}
where $\lambda_k=[\lambda_k^1;\lambda_k^2;\lambda_k^3;\lambda_k^4]\in\Re^{12}$, and $z_k\in\Re^{12}$ represents the infinitesimal variation of $(R_k,x_k,\Omega_k,v_k)$, given by $z_k=[\mathrm{logm}(R_k^T \delta R_k)^\vee;\delta x_k,\delta\Omega_k,\delta v_k]$. The matrix $A_k\in\Re^{12\times 12}$ is defined in terms of $(R_k,x_k,\Omega_k,v_k)$. Thus, necessary conditions for optimality are given by
\begin{align}
u^f_{k+1} &= -W_{f}^{-1}\lambda_{k}^2,\label{eqn:ufkp}\\
u^m_{k+1} &= -W_{m}^{-1}\lambda_{k}^4,\label{eqn:umkp}\\
\lambda_{k} &= A_{k+1}^T \lambda_{k+1}\label{eqn:updatelam}
\end{align}
together with the discrete equations of motion and the boundary conditions.

\textit{Computational approach:}
Necessary conditions for optimality are expressed in terms of a two point boundary problem. This problem is to find the optimal discrete flow, multipliers, and control inputs to satisfy the equations of motion, optimality conditions, multiplier equations, and boundary conditions simultaneously. We use a neighboring extremal method (see \citepri{Bry.BK75}). A nominal solution satisfying all of the necessary conditions except the boundary conditions is chosen. The unspecified initial multiplier is updated by successive linearization so as to satisfy the specified
terminal boundary conditions in the limit. This is also referred to as a shooting method. The main advantage of the neighboring extremal method is that the number of iteration variables is small.

The difficulty is that the extremal solutions are sensitive to small changes in the unspecified initial multiplier values. The nonlinearities also make it hard to construct an accurate estimate of sensitivity, thereby resulting in numerical ill-conditioning. Therefore, it is important to compute the sensitivities accurately to apply the neighboring extremal method. Here the optimality conditions
\refeqn{ufkp} and \refeqn{umkp} are substituted into the equations of motion and the multiplier equations, which are linearized to obtain sensitivity derivatives of an optimal solution with respect to boundary conditions.
Using this sensitivity, an initial guess of the unspecified initial conditions is iterated to satisfy the specified terminal conditions in the limit. Any type of Newton iteration can be applied. We use a line search with backtracking algorithm, referred to as Newton-Armijo iteration (see \citepri{Kel.BK95}).

\begin{figure}
\subfigure[Convergence rate]%
    {\includegraphics[width=0.32\textwidth]{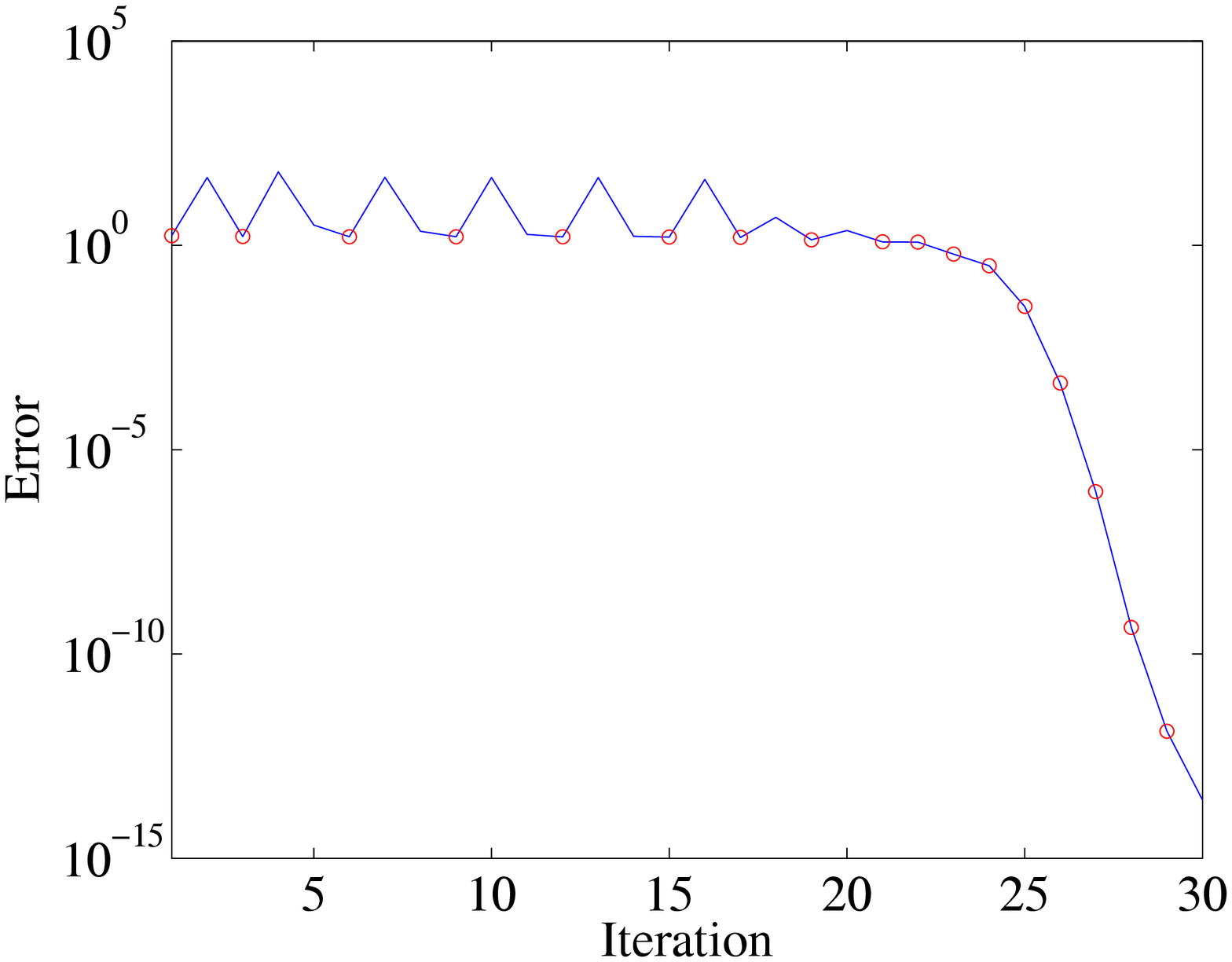}\label{fig:convopt}}
\hspace{0.002\textwidth}
\subfigure[Orbital radius change]%
    {\includegraphics[width=0.34\textwidth]{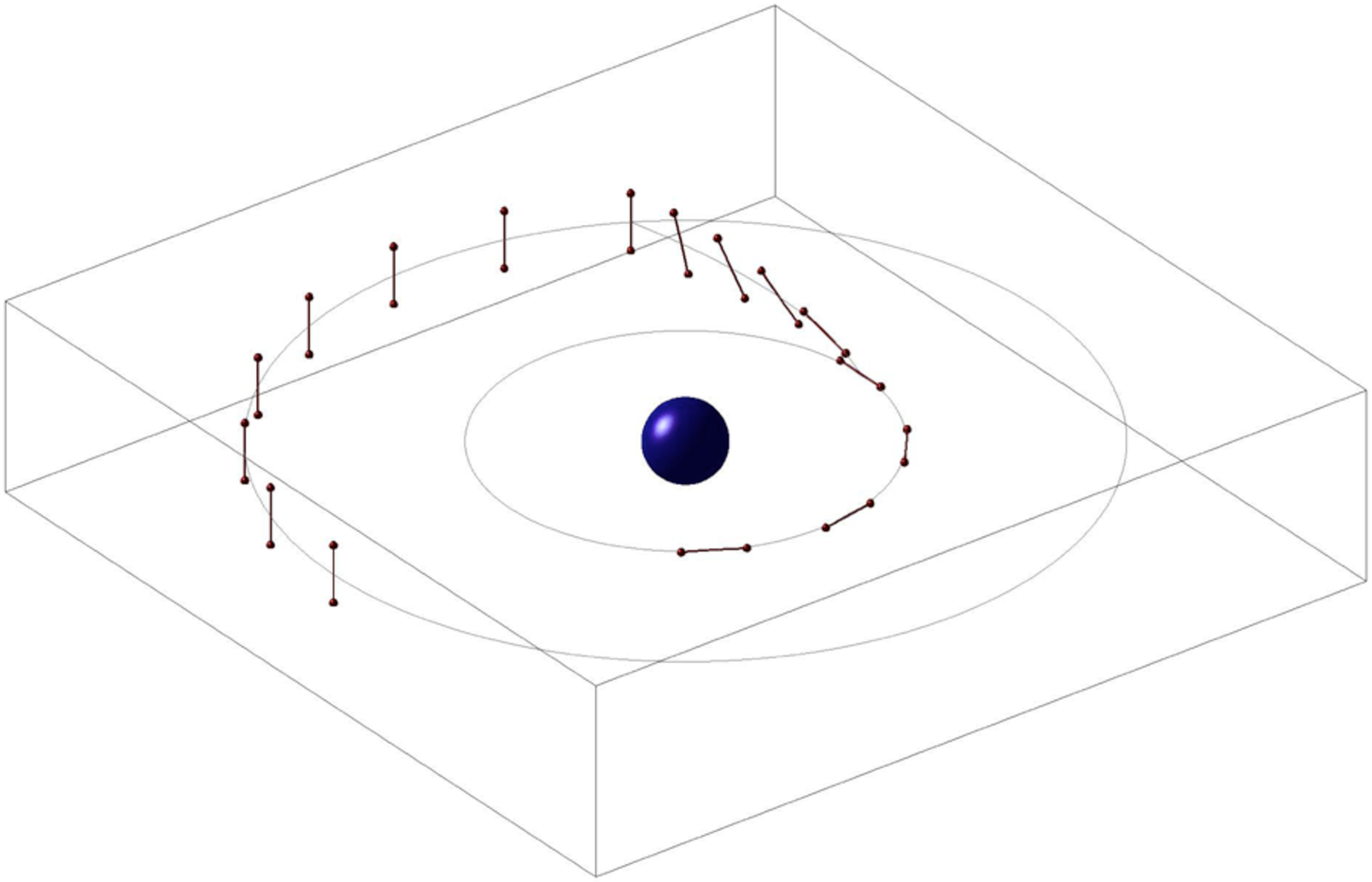}\label{fig:trajropt}}
\subfigure[Inclination angle change]%
    {\includegraphics[width=0.26\textwidth]{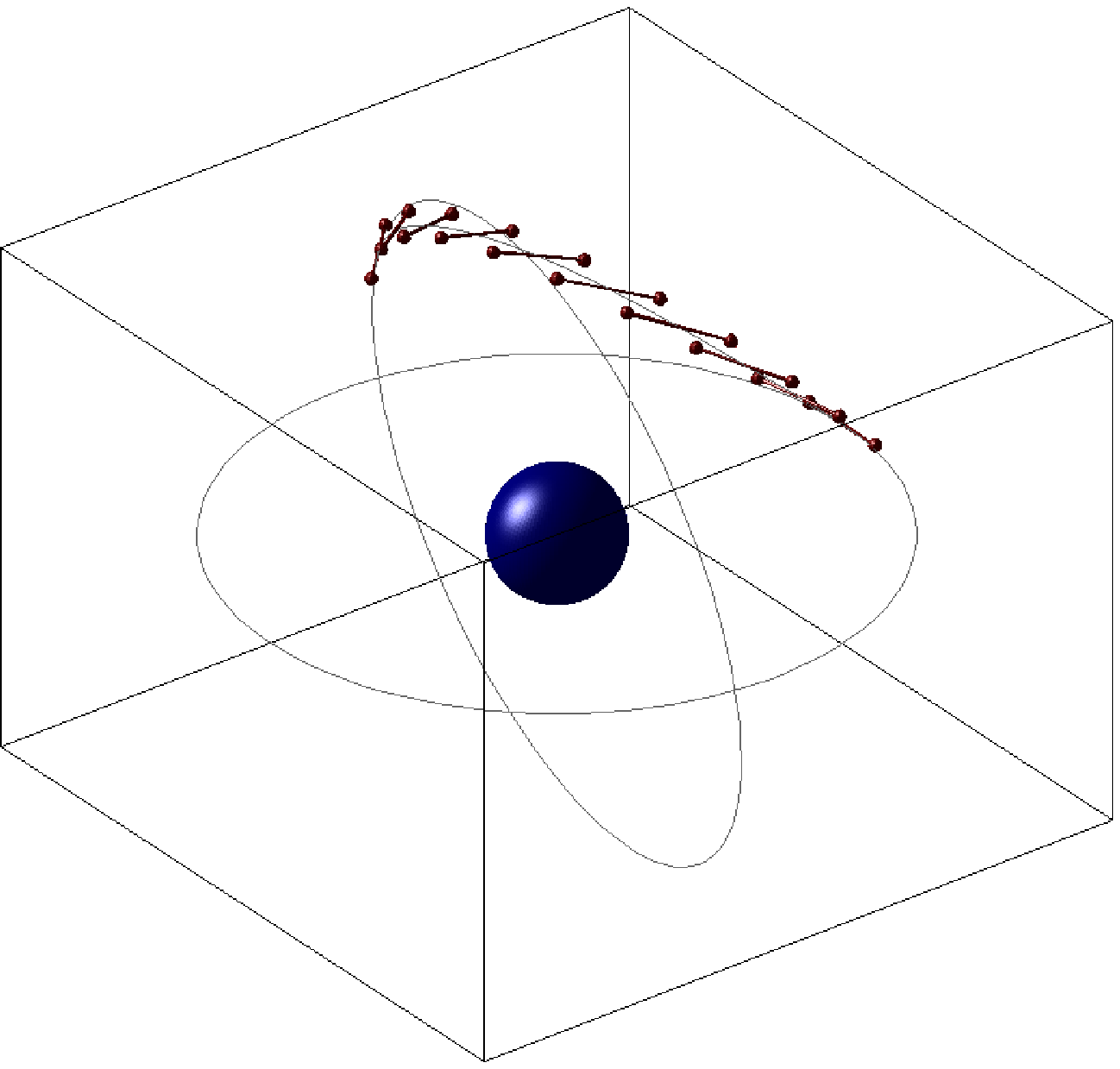}\label{fig:trajiopt}}
\caption{Optimal control of a rigid body}\label{fig:optrb}
\end{figure}

\reffig{optrb} shows optimized maneuvers, where a dumbbell spacecraft on a reference circular orbit is transferred to another circular orbits with different orbital radius and inclination angle. \reffig{convopt} shows the violation of the terminal boundary condition according to the number of iterations in a logarithmic scale. Red circles denote outer iterations in Newton-Armijo iteration to compute the sensitivity derivatives. The error in satisfaction of the terminal boundary condition converges quickly to machine precision after the solution is close to the local minimum at around 20th iteration. These convergence results are consistent with the quadratic convergence rates expected of Newton methods with accurately computed gradients.

The neighboring extremal method, also referred to as the shooting method, is numerically efficient in the sense that the number of optimization parameters is minimized. But, in general, this approach may be prone to numerical ill-conditioning (see ~\citepri{Bet.BK01}). A small change in the initial multiplier can cause highly nonlinear behavior of the terminal attitude and angular momentum. It is difficult to compute the gradient for Newton iterations accurately, and the numerical error may not converge. However, the numerical examples presented in this article show excellent numerical convergence properties. The dynamics of a rigid body arises from Hamiltonian mechanics, which have neutral stability, and its adjoint system is also neutrally stable. The proposed Lie group variational integrator and the discrete multiplier equations, obtained from variations expressed in the Lie algebra, preserve the neutral stability property numerically. Therefore the sensitivity derivatives are computed accurately.
\end{example}

\section{Controlled Lagrangian Method for Discrete Lagrangian Systems}
The method of controlled Lagrangians is a procedure for constructing feedback controllers for the stabilization of relative equilibria. It relies on choosing a parametrized family of controlled Lagrangians whose corresponding Euler--Lagrange flow are equivalent to the closed loop behavior of a Lagrangian system with external control forces. The condition that these two systems are equivalent result in matching conditions. Since the controlled system can now be viewed as a Lagrangian system with a modified Lagrangian, the global stability of the controlled system can be determined directly using Lyapunov stability analysis. 

This approach originated in \citepri{BLM1} and was then developed in \citepri{A,BLM2, BLM4, BLM3,BCLM,H1, H2}. A similar approach for Hamiltonian controlled systems was introduced and further studied in the work of Blankenstein, Ortega, van der Schaft, Maschke, Spong, and their collaborators (see, for example, \citepri{MaOrSc2000, OrSpGoBl2002} and related references). The two methods were shown to be equivalent in \citepri{ChBlLeMa2002} and a nonholonomic version was developed in \citepri{ZBM3,ZBM2002}, and \citepri{Bl2003}.

Since the method of controlled Lagrangians relies on relating the closed-loop dynamics of a controlled system with the Euler--Lagrange flow associated with a modified Lagrangian, it is natural to discretize this approach through the use of variational integrators. In \citepri{BloLeoMar.CDC05,BloLeoMar.CDC06}, a discrete theory of controlled Lagrangians was developed for variational integrators, and applied to the feedback stabilization of the unstable inverted equilibrium of the pendulum on a cart.

The pendulum on a cart is an example of an underactuated control problem, which has two degrees of freedom, given by the pendulum angle and the cart position. Only the cart position has control forces acting on it, and the stabilization of the pendulum has to be achieved indirectly through the coupling between the pendulum and the cart. The controlled Lagrangian is obtained by modifying the kinetic energy term, a process referred to as kinetic shaping. Similarly, it is possible to modify the potential energy term using potential shaping.

Since the pendulum on a cart model involves both a planar pendulum, and a cart that translates in one-dimension, the configuration space is a cylinder, $\mathbb{S}^1\times\mathbb{R}$.

\textit{Continuous kinetic shaping:}
The Lagrangian has the form kinetic minus potential energy
\begin{equation}\label{continuous_11_lagr.eqn}
L(q, \dot{q}) = \tfrac 12 \big[
	\alpha \dot \theta ^2 + 2 \beta (\theta) \dot \theta \dot s 
	+ \gamma \dot s ^2 
\big] - V (q),
\end{equation}
and the corresponding controlled Euler--Lagrange dynamics is
\begin{align}\label{EL_L_1.eqn} 
\frac{d}{dt} \frac{\partial L}{\partial \dot \theta} 
- \frac{\partial L}{\partial \theta} &= 0,
\\
\label{EL_L_2.eqn}
\frac{d}{dt} \frac{\partial L}{\partial \dot{s}} &= u,
\end{align} 
where $u$ is the control input. 

Since the potential energy is translationally invariant, \emph{i.e.,} $ V (q) = V (\theta) $, and the \emph{relative equilibria} $ \theta = \theta _e $, $ \dot{s} = \text{const} $ are unstable and given by non-degenerate critical points of $ V (\theta) $. To stabilize the relative equilibria $ \theta = \theta _e $, $ \dot{s} = \text{const} $ with respect to $\theta$, kinetic shaping is used. The controlled Lagrangian in this case is defined by
\begin{equation}\label{cont_controlled_lagrangian_11.eqn}
L^{\tau,\sigma}(q, \dot{q}) = 
L( \theta, \dot \theta, \dot{s} + \tau (\theta) \dot \theta) +
{\textstyle \frac12} \sigma \gamma (\tau (\theta) \dot \theta ) ^2, 
\end{equation} 
where $ \tau (\theta) = \kappa \beta (\theta) $. This velocity shift corresponds to a new choice of the horizontal space (see \citepri{BLM3} for details). The dynamics is just the Euler--Lagrange dynamics for controlled Lagrangian \eqref{cont_controlled_lagrangian_11.eqn},
\begin{align}\label{EL_CL_1.eqn}
\frac{d}{dt} \frac{\partial L ^{\tau, \sigma}}{\partial \dot \theta} 
- \frac{\partial L^{\tau, \sigma}}{\partial \theta} &= 0,
\\
\label{EL_CL_2.eqn}
\frac{d}{dt} \frac{\partial L^{\tau, \sigma}}{\partial \dot{s}} &= 0.
\end{align}
Lagrangian \eqref{cont_controlled_lagrangian_11.eqn} satisfies the simplified matching conditions of \citepri{BCLM} when the kinetic energy metric coefficient $\gamma$ in \eqref{continuous_11_lagr.eqn} is constant. 

Setting 
\(
u = - d \big( \gamma \tau (\theta) \dot{\theta} \big)/ dt 
\)
defines the control input, makes equations \eqref{EL_L_2.eqn} and \eqref{EL_CL_2.eqn} identical, and results in controlled momentum conservation by dynamics \eqref{EL_L_1.eqn} and \eqref{EL_L_2.eqn}. Setting $ \sigma = - 1/\gamma \kappa $ makes equations \eqref{EL_L_1.eqn} and \eqref{EL_CL_1.eqn} reduced on the controlled momentum level identical.

\textit{Discrete kinetic shaping:}
Here, we adopt the following notation: \begin{equation*}\label{notations.eqn}
q _{k+1/2} = \frac {q _k + q _{k+1}}{2}, \quad 
\Delta q _k = q _{k+1} - q _k, \quad q _k = (\theta _k, s _k).
\end{equation*} 
Then, a second-order accurate discrete Lagrangian is given by,
\[
L_d(q_k, q_{k+1})=hL(q_{k+1/2},\Delta q_k/h).
\]
The discrete dynamics is governed by the equations
\begin{align} 
\label{dcp1.eqn}
\pder{L _d (q _k,q  _{k+1})}{\theta _k} + 
\pder{L _d (q_{k-1},q  _k)}{\theta _k} &= 0,
\\
\label{dcp2.eqn}
\pder{L _d (q _k,q  _{k+1})}{s _k} + 
\pder{L _d (q _{k-1},q  _k)}{s _k} &= -u _k ,
\end{align}
where $ u _k $ is the control input. Similarly, the discrete controlled Lagrangian is,
\[
L_d^{\tau,\sigma}(q_k, q_{k+1})=hL^{\tau,\sigma}(q_{k+1/2},\Delta q_k/h),
\]
with discrete dynamics given by,
\begin{align} 
\label{controlled_cart_pendulum_1.eqn}
\pder{L _d  ^{ \tau , \sigma } (q _k,q  _{k+1})}{\theta _k} + 
\pder{L _d  ^{ \tau , \sigma } (q_{k-1},q  _k)}{\theta _k} &= 0,
\\
\label{controlled_cart_pendulum_2.eqn}
\pder{L _d  ^{ \tau , \sigma } (q _k,q  _{k+1})}{s _k} + 
\pder{L _d  ^{ \tau , \sigma } (q _{k-1},q  _k)}{s _k} &= 0. 
\end{align}
Equation \eqref{controlled_cart_pendulum_2.eqn} is equivalent to the
\emph{discrete controlled momentum conservation:}
\[
p _k = \mu ,
\] 
where 
\[
p _k = - \frac{\partial L _d ^{\tau,\sigma} (q _k, q _{k+1})}{\partial s _k} 
=  \frac{
	(1 + \gamma \kappa) \beta (\theta _{k+1/2} ) 
	\Delta \theta _k + \gamma \Delta s _k} {h}.
\] 

Setting
\[
u _k = - \frac{\gamma \Delta \theta _k \tau ( \theta _{k+1/2} ) 
- \gamma \Delta \theta _{k-1} \tau ( \theta _{k-1/2})}{h }
\] 
makes equations \eqref{dcp2.eqn} and \eqref{controlled_cart_pendulum_2.eqn}
identical and allows 
one to represent the discrete momentum equation
\eqref{dcp2.eqn} as the discrete momentum conservation law 
$p _k = p.$

The condition that (\ref{dcp1.eqn}-\ref{dcp2.eqn}) are equivalent to (\ref{controlled_cart_pendulum_1.eqn}-\ref{controlled_cart_pendulum_2.eqn}) yield the discrete matching conditions. The dynamics determined by equations (\ref{dcp1.eqn}-\ref{dcp2.eqn}) restricted to the momentum level $ p _k = p $ is equivalent to the dynamics of equations
(\ref{controlled_cart_pendulum_1.eqn}-\ref{controlled_cart_pendulum_2.eqn}) restricted to the momentum level $
p _k = \mu $ if and only if the matching conditions 
\[
\sigma 
= - \frac{1}{\gamma \kappa}, \qquad 
\mu = \frac{p}{1 + \gamma \kappa},
\] hold.

\begin{figure}[ht]
\hspace{-1.15em}
\begin{overpic}
[scale=.51]
{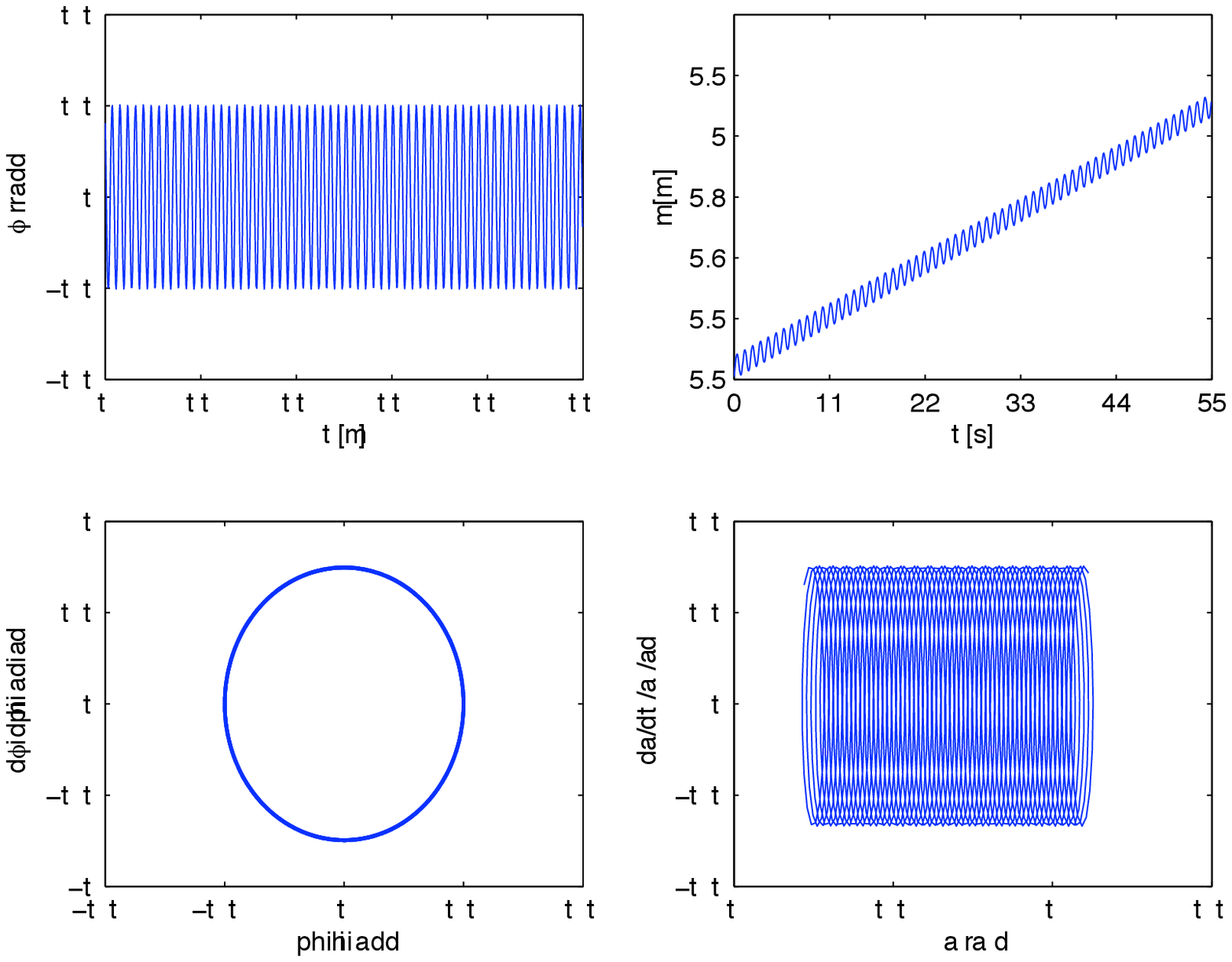}
\put(28.3,10.55)
	{\begin{rotate}{0}
	{\includegraphics[width=.015\textwidth]
	{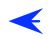}}
	\end{rotate}
	}
\put(31.5,31.67)
	{\begin{rotate}{180}
	{\includegraphics[width=.015\textwidth]
	{dcs_fig7a}}
	\end{rotate}
	}
\end{overpic} 
\caption{\label{fig:discrete_kinetic_nodiss}Discrete controlled dynamics with kinetic shaping and without dissipation. The discrete controlled system stabilizes the $\theta$ motion about the equilibrium, but the $s$ dynamics is not stabilized; since there is no dissipation, the oscillations are sustained.}
\end{figure}

\textit{Numerical example:} Simulating the behavior of the discrete controlled Lagrangian system
involves viewing equations (\ref{controlled_cart_pendulum_1.eqn}-\ref{controlled_cart_pendulum_2.eqn}) as an implict update map $(q_{k-2},
q_{k-1})\mapsto(q_{k-1},q_k)$. This presupposes that the initial
conditions are given in the form $(q_0,q_1)$; however it is generally
preferable to specify the initial conditions as $(q_0,\dot q_0)$. 
This is achieved by solving the boundary condition,
\[
\frac{\partial L}{\partial \dot q}(q_0,\dot q_0) + D_1 L_d(q_0,q_1) +
F_d^{-}(q_0,q_1)
= 0,
\]
for $q_1$. Once the initial conditions are expressed in the form $(q_0,q_1)$, the discrete evolution can be obtained using the implicit update
map.

We first consider the case of kinetic shaping on a level surface, when $\kappa$ is twice the critical value, and without dissipation. Here, $h=0.05\,\textrm{sec}$, $m=0.14\,\textrm{kg}$, $M=0.44\,\textrm{kg}$, and $l=0.215\,\textrm{m}$. As shown in Figure~\ref{fig:discrete_kinetic_nodiss}, the $\theta$ dynamics is stabilized, but since there is no dissipation, the oscillations are sustained. The $s$ dynamics exhibits both a drift and oscillations, as potential shaping is necessary to stabilize the translational dynamics.

\section{Future Directions}
\textit{Discrete Receding Horizon Optimal Control:} The existing work on discrete optimal control has been primarily focused on constructing the optimal trajectory in an open loop sense. In practice, model uncertainty and actuation errors necessitate the use of feedback control, and it would be interesting to extend the existing work on optimal control of discrete systems to the feedback setting by adopting a receding horizon approach.

\textit{Discrete State Estimation:}
In feedback control, one typically assumes complete knowledge regarding the state of the system, an assumption that is often unrealistic in practice. The general problem of state estimation in the context of discrete mechanics would rely on good numerical methods for quantifying the propagation of uncertainty by solving the Liouville equation. In the setting of Hamiltonian systems, the solution of the Liouville equation can be solved by the method of characteristics (\citepri{ScHsPaViMa2007}). This implies that a collocational approach (\citepri{Xi2007}) combined with Lie group variational integrators, and interpolation based on noncommutative harmonic analysis on Lie groups could yield an efficient means of propagating uncertainty, and serve as the basis of a discrete state estimation algorithm.

\textit{Forced Symplectic-Energy-Momentum Variational Integrators:}
One of the motivations for studying the control of Lagrangian systems using the method of controlled Lagrangians is that the method provides a natural candidate Lyapunov function to study the global stability properties of the controlled system. In the discrete theory, this approach is complicated by the fact that the energy of a discrete Lagrangian system is not exactly conserved, but rather oscillates in a bounded fashion.

This can be addressed by considering the symplectic-energy-momentum (\citepri{KaMaOr1999}) analogue of the discrete Lagrange-d'Alembert principle,
\[
\delta\sum_{k=0}^{N_1}L_d(q_k,q_{k+1},h_k)=\sum_{k=0}^{N-1} [F_d^{-}(q_k,q_{k+1},h_k)\cdot\delta q_k+F_d^+(q_k,q_{k+1},h_k)\cdot\delta q_{k+1}],
\]
where the timestep $h_k$ is allowed to vary, and is chosen to satisfy the variational principle. The variations in $h_k$ yield an Euler--Lagrange equation that reduces to the conservation of discrete energy in the absence of external forces. By developing a theory of controlled Lagrangians around a geometric integrator based on the symplectic-energy-momentum version of the Lagrange-d'Alembert principle, one would potentially be able to use Lyapunov techniques to study the global stability of the resulting numerical control algorithms.

\nociterev{Bl2003,BuLe2005,HaLuWa2006,IsMuNoZa2000,LeRe2004,MarRat.BK99,MaWe2001,San.AN92}
\bibliographystylepri{plainnat}
\bibliographypri{dcs}
\bibliographystylerev{plainnat}
\bibliographyrev{dcs}

\begin{thebibliography}{40}
\providecommand{\natexlab}[1]{#1}
\providecommand{\url}[1]{\texttt{#1}}
\expandafter\ifx\csname urlstyle\endcsname\relax
  \providecommand{\doi}[1]{doi: #1}\else
  \providecommand{\doi}{doi: \begingroup \urlstyle{rm}\Url}\fi

\bibitem[Auckly et~al.(2000)Auckly, Kapitanski, and White]{A}
D.~Auckly, L.~Kapitanski, and W.~White.
\newblock Control of nonlinear underactuated systems.
\newblock \emph{Commun. Pure Appl. Math.}, 53:\penalty0 354--369, 2000.

\bibitem[Betts(2001)]{Bet.BK01}
J.~T. Betts.
\newblock \emph{Practical Methods for Optimal Control Using Nonlinear
  Programming}.
\newblock SIAM, 2001.

\bibitem[Bloch(2003)]{Bl2003}
{A. M.} Bloch.
\newblock \emph{Nonholonomic Mechanics and Control}, volume~24 of
  \emph{Interdisciplinary Appl. Math.}
\newblock Springer-Verlag, 2003.

\bibitem[Bloch et~al.(1997)Bloch, Leonard, and Marsden]{BLM1}
{A. M.} Bloch, N.~Leonard, and {J. E.} Marsden.
\newblock Matching and stabilization using controlled {L}agrangians.
\newblock In \emph{Proceedings of the IEEE Conference on Decision and Control},
  pages 2356--2361, 1997.

\bibitem[Bloch et~al.(1998)Bloch, Leonard, and Marsden]{BLM2}
{A. M.} Bloch, N.~Leonard, and {J. E.} Marsden.
\newblock Matching and stabilization by the method of controlled {L}agrangians.
\newblock In \emph{Proceedings of the IEEE Conference on Decision and Control},
  pages 1446--1451, 1998.

\bibitem[Bloch et~al.(1999)Bloch, Leonard, and Marsden]{BLM4}
{A. M.} Bloch, N.~Leonard, and {J. E.} Marsden.
\newblock Potential shaping and the method of controlled {L}agrangians.
\newblock In \emph{Proceedings of the IEEE Conference on Decision and Control},
  pages 1652--1657, 1999.

\bibitem[Bloch et~al.(2000)Bloch, Leonard, and Marsden]{BLM3}
{A. M.} Bloch, {N. E.} Leonard, and {J. E.} Marsden.
\newblock Controlled {L}agrangians and the stabilization of mechanical systems
  {I}: The first matching theorem.
\newblock \emph{IEEE Trans. on Systems and Control}, 45:\penalty0 2253--2270,
  2000.

\bibitem[Bloch et~al.(2001)Bloch, Chang, Leonard, and Marsden]{BCLM}
{A. M.} Bloch, {D. E.} Chang, {N. E.} Leonard, and {J. E.} Marsden.
\newblock Controlled {L}agrangians and the stabilization of mechanical systems
  {II}: Potential shaping.
\newblock \emph{IEEE Trans. on Autom. Contr.}, 46:\penalty0 1556--1571, 2001.

\bibitem[Bloch et~al.(2005)Bloch, Leok, Marsden, and Zenkov]{BloLeoMar.CDC05}
A.~M. Bloch, M.~Leok, J.~E. Marsden, and D.~V. Zenkov.
\newblock Controlled {L}agrangians and stabilization of the discrete
  cart-pendulum system.
\newblock In \emph{Proceedings of the IEEE Conference on Decision and Control},
  pages 6579--6584, 2005.

\bibitem[Bloch et~al.(2006)Bloch, Leok, Marsden, and Zenkov]{BloLeoMar.CDC06}
A.~M. Bloch, M.~Leok, J.~E. Marsden, and D.~V. Zenkov.
\newblock Controlled {L}agrangians and potential shaping for stabilization of
  discrete mechanical systems.
\newblock In \emph{Proceedings of the IEEE Conference on Decision and Control},
  pages 3333--3338, 2006.

\bibitem[Bryson and Ho(1975)]{Bry.BK75}
A.~E. Bryson and Y.~Ho.
\newblock \emph{Applied Optimal Control}.
\newblock Hemisphere Publishing Corporation, 1975.

\bibitem[Chang et~al.(2002)Chang, Bloch, Leonard, Marsden, and
  Woolsey]{ChBlLeMa2002}
{D-E.} Chang, {A. M.} Bloch, {N. E.} Leonard, {J. E.} Marsden, and C.~Woolsey.
\newblock The equivalence of controlled {L}agrangian and controlled
  {H}amiltonian systems.
\newblock \emph{Control and the Calculus of Variations (special issue dedicated
  to {J. L.} Lions)}, 8:\penalty0 393--422, 2002.

\bibitem[Hairer et~al.(2006)Hairer, Lubich, and Wanner]{HaLuWa2006}
E.~Hairer, C.~Lubich, and G.~Wanner.
\newblock \emph{Geometric Numerical Integration}, volume~31 of \emph{Springer
  Series in Computational Mathematics}.
\newblock Springer-Verlag, second edition, 2006.

\bibitem[Hamberg(1999)]{H1}
J.~Hamberg.
\newblock General matching conditions in the theory of controlled
  {L}agrangians.
\newblock In \emph{Proceedings of the IEEE Conference on Decision and Control},
  pages 2519--2523, 1999.

\bibitem[Hamberg(2000)]{H2}
J.~Hamberg.
\newblock Controlled {L}agrangians, symmetries and conditions for strong
  matching.
\newblock In \emph{Lagrangian and Hamiltonian Methods for Nonlinear Control}.
  Elsevier, 2000.

\bibitem[Hussein et~al.(2006)Hussein, Leok, Sanyal, and Bloch]{HuLeSaBl2006}
{I. I.} Hussein, M.~Leok, {A. K.} Sanyal, and {A.M.} Bloch.
\newblock A discrete variational integrator for optimal control problems in
  {$SO(3)$}.
\newblock In \emph{Proceedings of the IEEE Conference on Decision and Control},
  pages 6636--6641, 2006.

\bibitem[Iserles et~al.(2000)Iserles, Munthe-Kaas, {N\o rsett}, and
  Zanna]{IsMuNoZa2000}
A.~Iserles, H.~Munthe-Kaas, {S. P.} {N\o rsett}, and A.~Zanna.
\newblock Lie-group methods.
\newblock In \emph{Acta Numerica}, volume~9, pages 215--365. Cambridge
  University Press, 2000.

\bibitem[Junge et~al.(2005)Junge, Marsden, and Ober-Bl{\"o}baum]{JuMaOb2005}
O.~Junge, J.~E. Marsden, and S.~Ober-Bl{\"o}baum.
\newblock Discrete mechanics and optimal control.
\newblock In \emph{IFAC Congress}, Praha, 2005.

\bibitem[Kane et~al.(1999)Kane, Marsden, and Ortiz]{KaMaOr1999}
C.~Kane, J.~E. Marsden, and M.~Ortiz.
\newblock Symplectic-energy-momentum preserving variational integrators.
\newblock \emph{J. Math. Phys.}, 40\penalty0 (7):\penalty0 3353--3371, 1999.

\bibitem[Kane et~al.(2000)Kane, Marsden, Ortiz, and West]{KaMaOrWe2000}
C.~Kane, {J. E.} Marsden, M.~Ortiz, and M.~West.
\newblock Variational integrators and the {N}ewmark algorithm for conservative
  and dissipative mechanical systems.
\newblock \emph{Int. J. Numer. Meth. Eng.}, 49\penalty0 (10):\penalty0
  1295--1325, 2000.

\bibitem[Kelley(1995)]{Kel.BK95}
C.~T. Kelley.
\newblock \emph{Iterative Methods for Linear and Nonlinear Equations}.
\newblock SIAM, 1995.

\bibitem[Lee et~al.(2005{\natexlab{a}})Lee, Leok, and McClamroch]{ACC06}
T.~Lee, M.~Leok, and N.~H. McClamroch.
\newblock Attitude maneuvers of a rigid spacecraft in a circular orbit.
\newblock In \emph{Proceedings of the American Control Conference}, pages
  1742--1747, 2005{\natexlab{a}}.

\bibitem[Lee et~al.(2005{\natexlab{b}})Lee, Leok, and McClamroch]{CCA05}
T.~Lee, M.~Leok, and {N. H.} McClamroch.
\newblock A {L}ie group variational integrator for the attitude dynamics of a
  rigid body with applications to the 3{D} pendulum.
\newblock In \emph{Proceedings of the IEEE Conference on Control Applications},
  pages 962--967, 2005{\natexlab{b}}.

\bibitem[Lee et~al.(2006)Lee, Leok, and McClamroch]{CDC06.opt}
T.~Lee, M.~Leok, and N.~H. McClamroch.
\newblock Optimal control of a rigid body using geometrically exact
  computations on {S}{E}(3).
\newblock In \emph{Proceedings of the IEEE Conference on Decision and Control},
  pages 2710--2715, 2006.

\bibitem[Lee et~al.(2007{\natexlab{a}})Lee, Leok, and McClamroch]{CMA07}
T.~Lee, M.~Leok, and N.~H. McClamroch.
\newblock Lie group variational integrators for the full body problem.
\newblock \emph{Computer Methods in Applied Mechanics and Engineering},
  196:\penalty0 2907--2924, May 2007{\natexlab{a}}.

\bibitem[Lee et~al.(2007{\natexlab{b}})Lee, Leok, and McClamroch]{CMDA07}
T.~Lee, M.~Leok, and N.~H. McClamroch.
\newblock Lie group variational integrators for the full body problem in
  orbital mechanics.
\newblock \emph{Celestial Mechanics and Dynamical Astronomy}, 98\penalty0
  (2):\penalty0 121--144, June 2007{\natexlab{b}}.

\bibitem[Leimkuhler and Reich(2004)]{LeRe2004}
B.~Leimkuhler and S.~Reich.
\newblock \emph{Simulating {H}amiltonian Dynamics}, volume~14 of
  \emph{Cambridge Monographs on Applied and Computational Mathematics}.
\newblock Cambridge University Press, 2004.

\bibitem[Leok(2004)]{Leo.Phd04}
M.~Leok.
\newblock \emph{Foundations of {C}omputational {G}eometric {M}echanics}.
\newblock PhD thesis, California {I}nstittute of {T}echnology, 2004.

\bibitem[Marsden and Ratiu(1999)]{MarRat.BK99}
J.~E. Marsden and T.~S. Ratiu.
\newblock \emph{Introduction to Mechanics and Symmetry}, volume~17 of
  \emph{Texts in Applied Mathematics}.
\newblock Springer-Verlag, second edition, 1999.

\bibitem[Marsden and West(2001)]{MaWe2001}
{J. E.} Marsden and M.~West.
\newblock Discrete mechanics and variational integrators.
\newblock In \emph{Acta Numerica}, volume~10, pages 317--514. Cambridge
  University Press, 2001.

\bibitem[Marsden et~al.(1999)Marsden, Pekarsky, and Shkoller]{MaPeSh1999}
{J. E.} Marsden, S.~Pekarsky, and S.~Shkoller.
\newblock Discrete {E}uler--{P}oincar\'e and {L}ie--{P}oisson equations.
\newblock \emph{Nonlinearity}, 12\penalty0 (6):\penalty0 1647--1662, 1999.

\bibitem[Maschke et~al.(2001)Maschke, Ortega, and {van der Schaft}]{MaOrSc2000}
B.~Maschke, R.~Ortega, and A.~{van der Schaft}.
\newblock Energy-based {L}yapunov functions for forced {H}amiltonian systems
  with dissipation.
\newblock \emph{IEEE Trans. on Autom. Contr.}, 45:\penalty0 1498--1502, 2001.

\bibitem[Moser and Veselov(1991)]{MosVes.CMP91}
J.~Moser and A.~P. Veselov.
\newblock Discrete versions of some classical integrable systems and
  factorization of matrix polynomials.
\newblock \emph{Communications in Mathematical Physics}, 139:\penalty0
  217--243, 1991.

\bibitem[Ortega et~al.(2002)Ortega, Spong, {G\'{o}mez-Estern}, and
  Blankenstein]{OrSpGoBl2002}
R.~Ortega, {M. W.} Spong, F.~{G\'{o}mez-Estern}, and G.~Blankenstein.
\newblock Stabilization of a class of underactuated mechanical systems via
  interconnection and damping assignment.
\newblock \emph{IEEE Trans. on Autom. Contr.}, 47:\penalty0 1218--1233, 2002.

\bibitem[Sanz-Serna(1992)]{San.AN92}
J.~M. Sanz-Serna.
\newblock Symplectic integrators for {H}amiltonian problems: an overview.
\newblock In \emph{Acta Numerica}, volume~1, pages 243--286. Cambridge
  University Press, 1992.

\bibitem[Scheeres et~al.(2006)Scheeres, Fahnestock, Ostro, Margot, Benner,
  Broschart, Bellerose, Giorgini, Nolan, Magri, Pravec, Scheirich, Rose,
  Jurgens, Jong, and Suzuki]{Sch.Sci06}
D.~J. Scheeres, E.~G. Fahnestock, S.~J. Ostro, J.~L. Margot, L.~A.~M. Benner,
  S.~B. Broschart, J.~Bellerose, J.~D. Giorgini, M.~C. Nolan, C.~Magri,
  P.~Pravec, P.~Scheirich, R.~Rose, R.~F. Jurgens, E.~M.~De Jong, and
  S.~Suzuki.
\newblock Dynamical configuration of binary near-{E}arth asteroid (66391) 1999
  {K}{W}4.
\newblock \emph{Science}, 314:\penalty0 1280--1283, 2006.

\bibitem[Scheeres et~al.(2007)Scheeres, Hsiao, Park, Villac, and
  Maruskin]{ScHsPaViMa2007}
{D. J.} Scheeres, {F.-Y.} Hsiao, {R. S.} Park, {B. F.} Villac, and {J. M.}
  Maruskin.
\newblock Fundamental limits on spacecraft orbit uncertainty and distribution
  propagation.
\newblock \emph{Journal of the Astronautical Sciences}, 2007.
\newblock to appear.

\bibitem[Xiu(2007)]{Xi2007}
D.~Xiu.
\newblock Efficient collocational approach for parametric uncertainty analysis.
\newblock \emph{Comm. Comput. Phys.}, 2:\penalty0 293--309, 2007.

\bibitem[Zenkov et~al.(2000)Zenkov, Bloch, Leonard, and {J. E.}Marsden]{ZBM3}
{D. V.} Zenkov, {A. M.} Bloch, {N. E.} Leonard, and {J. E.}Marsden.
\newblock Matching and stabilization of low-dimensional nonholonomic systems.
\newblock In \emph{Proceedings of the IEEE Conference on Decision and Control},
  pages 1289--1295, 2000.

\bibitem[Zenkov et~al.(2002)Zenkov, Bloch, Leonard, and {J.
  E.}Marsden]{ZBM2002}
{D. V.} Zenkov, {A. M.} Bloch, {N. E.} Leonard, and {J. E.}Marsden.
\newblock Flat nonholonomic matching.
\newblock In \emph{Proceedings of the American Control Conference}, pages
  2812--2817, 2002.

\end{thebibliography}


\begin{thebibliography}{8}
\providecommand{\natexlab}[1]{#1}
\providecommand{\url}[1]{\texttt{#1}}
\expandafter\ifx\csname urlstyle\endcsname\relax
  \providecommand{\doi}[1]{doi: #1}\else
  \providecommand{\doi}{doi: \begingroup \urlstyle{rm}\Url}\fi

\bibitem[Bloch(2003)]{Bl2003}
{A. M.} Bloch.
\newblock \emph{Nonholonomic Mechanics and Control}, volume~24 of
  \emph{Interdisciplinary Appl. Math.}
\newblock Springer-Verlag, 2003.

\bibitem[Bullo and Lewis(2005)]{BuLe2005}
F.~Bullo and {A. D.} Lewis.
\newblock \emph{Geometric control of mechanical systems}, volume~49 of
  \emph{Texts in Applied Mathematics}.
\newblock Springer-Verlag, New York, 2005.
\newblock Modeling, analysis, and design for simple mechanical control systems.

\bibitem[Hairer et~al.(2006)Hairer, Lubich, and Wanner]{HaLuWa2006}
E.~Hairer, C.~Lubich, and G.~Wanner.
\newblock \emph{Geometric Numerical Integration}, volume~31 of \emph{Springer
  Series in Computational Mathematics}.
\newblock Springer-Verlag, second edition, 2006.

\bibitem[Iserles et~al.(2000)Iserles, Munthe-Kaas, {N\o rsett}, and
  Zanna]{IsMuNoZa2000}
A.~Iserles, H.~Munthe-Kaas, {S. P.} {N\o rsett}, and A.~Zanna.
\newblock Lie-group methods.
\newblock In \emph{Acta Numerica}, volume~9, pages 215--365. Cambridge
  University Press, 2000.

\bibitem[Leimkuhler and Reich(2004)]{LeRe2004}
B.~Leimkuhler and S.~Reich.
\newblock \emph{Simulating {H}amiltonian Dynamics}, volume~14 of
  \emph{Cambridge Monographs on Applied and Computational Mathematics}.
\newblock Cambridge University Press, 2004.

\bibitem[Marsden and Ratiu(1999)]{MarRat.BK99}
J.~E. Marsden and T.~S. Ratiu.
\newblock \emph{Introduction to Mechanics and Symmetry}, volume~17 of
  \emph{Texts in Applied Mathematics}.
\newblock Springer-Verlag, second edition, 1999.

\bibitem[Marsden and West(2001)]{MaWe2001}
{J. E.} Marsden and M.~West.
\newblock Discrete mechanics and variational integrators.
\newblock In \emph{Acta Numerica}, volume~10, pages 317--514. Cambridge
  University Press, 2001.

\bibitem[Sanz-Serna(1992)]{San.AN92}
J.~M. Sanz-Serna.
\newblock Symplectic integrators for {H}amiltonian problems: an overview.
\newblock In \emph{Acta Numerica}, volume~1, pages 243--286. Cambridge
  University Press, 1992.

\end{thebibliography}
\end{document}